\documentclass[10pt]{amsart}
\usepackage{amsfonts}
\usepackage{graphicx}
\usepackage{amscd}
\usepackage{amsmath,amsfonts,amssymb,amsthm,mathrsfs,xypic}
\xyoption{curve}
\usepackage{fancybox}
\newtheorem{thm}{Theorem}[section]

\newtheorem{cor}[thm]{Corollary}

\newtheorem{lem}[thm]{Lemma}

\newtheorem{exm}[thm]{Example}

\newcommand{\D}{{\rm D}}

\newtheorem{prop}[thm]{Proposition}
\newtheorem{rem}[thm]{Remark}

\newtheorem{defn}[thm]{Definition}
\numberwithin{equation}{section}
\newcommand{\m}{\Lambda}
\newcommand{\s}{\hfill\blacksquare}
\newcommand{\add}{\operatorname{add}}

\newcommand{\cok}{\operatorname{Coker}}
\newcommand{\Hom}{\operatorname{Hom}}
\newcommand{\Ext}{\operatorname{Ext}}

\begin{document}
\title [Monic modules and semi-Gorenstein-projective modules]{Monic modules and semi-Gorenstein-projective modules}
\author [Pu Zhang] {Pu Zhang  \\ School of Mathematical Sciences, Shanghai Jiao Tong University \\
Shanghai 200240, P. R. China}
\thanks{Supported by National Natural Science Foundation of China, Grant No. 12131015, 11971304.}
\thanks{pzhang$\symbol{64}$sjtu.edu.cn}
\maketitle
\vspace {-20pt}
\begin{abstract} The category ${\rm gp}(\m)$ of Gorenstein-projective modules over tensor algebra $\m = A\otimes_kB$
can be described as the monomorphism category ${\rm mon}(B, {\rm gp}(A))$ of $B$ over ${\rm gp}(A)$.
In particular, Gorenstein-projective $\m$-modules are monic. In this paper,
we find the similar relation between semi-Gorenstein-projective $\m$-modules
and $A$-modules, via monic modules, namely, ${\rm mon}(B, \ ^\perp A) = {\rm mon}(B, A)\cap \ ^\perp \m.$
Using this, it is proved that if $A$ is weakly Gorenstein, then $\m$ is weakly Gorenstein if and only each semi-Gorenstein-projective $\m$-modules are monic;
and that if $B = kQ$ with $Q$ a finite acyclic quiver, then $\m$ is weakly Gorenstein if and only if $A$ is weakly Gorenstein.
However, this relation itself does not answer the question whether there exist double semi-Gorenstein-projective $\m$-modules which are not monic. Using the recent discovered examples of double semi-Gorenstein-projective $A$-modules which are not torsionless,
we positively answer this question, by explicitly
constructing a class of double semi-Gorenstein-projective $T_2(A)$-modules with one parameter such that they are not monic, and hence not torsionless.
The corresponding results are obtained also for the monic modules and semi-Gorenstein-projective modules over the triangular matrix algebras given by bimodules.

\vskip5pt


\vskip5pt

 \end{abstract}

\vskip10pt

\section{\bf Introduction}

Monic modules, defined on tensor products $\m = A\otimes_kB$, or on matrix algebras
$\m = \left(\begin{smallmatrix}A&M\\0&B\end{smallmatrix}\right)$ of bimodule modules $M$, built a bridge
between Gorenstein-projective $\m$-modules and Gorenstein-projective $A$-modules (Theorems \ref{gp} and \ref{gpovermatrixalg}).
In particular, Gorenstein-projective $\m$-modules are monic, in the both cases.
This paper is to show that they also play an important role in the study of semi-Gorenstein-projective modules. In the both cases, we will give sufficient and necessary conditions such that $\m$ is weakly Gorenstein,
and positively answer the question whether there exist double semi-Gorenstein-projective $\m$-modules which are not monic, and hence not torsionless, and not Goresnstein-projective.

\subsection{} Let $A$ be an Artin algebra. All the modules in this paper are finitely generated, and we start from left modules.
Let $A\mbox{-}{\rm mod}$ be the category of left $A$-modules. For $M\in A\mbox{-}{\rm mod}$, denote by
$\add(M)$ the full subcategory of $A\mbox{-}{\rm mod}$ of direct summands of a direct sum of copies of $M$;
by ${}^\perp M$ the full subcategory of $A\mbox{-}{\rm mod}$ of modules $X$ with $\Ext_A^i(X, M) = 0$
for $i\ge 1;$ and by $M^\perp$ the full subcategory of $A\mbox{-}{\rm mod}$ of modules $X$ with $\Ext_A^i(M, X) = 0$
for $i\ge 1.$

\vskip5pt

Let $M^*$ denote {\it the $A$-dual} $\Hom_A(M, A)$ of $M$.
Denote by $\phi_M: M \longrightarrow M^{**}$ the canonical $A$-map, defined by $\phi_M(m)(f) = f(m)$ for $m\in M$ and $f\in M^*$.
A module $M$ is {\it torsionless} if it is a submodule of a projective module, or, equivalently, $\phi_M$ is a monomorphism; and $M$ is {\it reflexive} if $\phi_M$ is an isomorphism.

\vskip5pt

A module $M$ is {\it semi-Gorenstein-projective} if $M\in {}^\perp A$; and $M$ will be called {\it double semi-Gorenstein-projective}, if both $M$ and $M^*$ are semi-Gorenstein-projective.
By definition, {\it a Gorenstein-projective} module is double semi-Gorenstein-projective and reflexive.
This is introduced by Auslander and Bridger [AB], under the name of {\it modules of $G$-dimension zero}, and it is equivalent to
the definition in terms of {\it complete projective resolution} given by Enochs and Jenda ([EJ1], [EJ2]).
For the equivalence we refer to [AM, p.398] (where it is called {\it a total reflexive module}) and  [Chr, Theorem 4.2.6]. Denote by ${\rm gp}(A)$ the full subcategory of $A\mbox{-}{\rm mod}$ of Gorenstein-projective modules. Thus
$\add(A) \subseteq {\rm gp}(A) \subseteq {}^\perp A.$

\subsection{} Avramov and Martsinkovsky [AM, p.398] has proposed the independence problem of the total reflexivity.
In fact, the known examples of semi-Gorenstein-projective modules which are not Gorenstein-projective are few and complicated.
The first examples of reflexive semi-Gorenstein-projective modules which are not Gorenstein-projective, and the first examples of reflexive modules $M$ with $M^*$ semi-Gorenstein-projective
such that $M$ are not semi-Gorenstein-projective, are discovered by Jorgensen and \c{S}ega [JS]; and the first examples of double
semi-Gorenstein-projective modules which are not torsionless, are recently founded in [RZ2, RZ3]. Putting together,
this solves the independence problem of the total reflexivity. Note that the first examples of semi-Gorenstein-projective modules which are not
Gorenstein-projective over noncommutative algebras, are presented by Marczinzik [M2].

\subsection{} Let $A$ and $B$ be finite-dimensional algebras over field $k$. Since Cartan-Eilenberg [CE], modules over tensor algebra $\m = A\otimes_k B$ have got interest.
They are complicated in the sense that $\m$-mod can not be controlled by $U\otimes_kV$ with $U\in A$-mod and $V\in B$-mod. However,
if $B$ is given by a bound quiver $(Q, I)$, one can study $\m$-module by taking the advantage of the representations of quivers over algebra $A$
([RS1-RS3], [S2-S5], [KLM1, KLM2], [LZ1, LZ2], [RZ1], [ZX]), i.e.,  any $\m$-module can be identified with {\it a
representation $(X_i, \ X_{\alpha}, \ i\in Q_0, \ \alpha\in Q_1)$ of $(Q, I)$ over $A$}, where each $X_i\in A$-mod, and each $X_\alpha: X_{s(\alpha)}\longrightarrow X_{e(\alpha)}$ is an
$A$-map, such that $X_\alpha$'s satisfy all the relations which generate $I$.

\vskip5pt

When $Q$ is finite acyclic and $I$ is generated by monomial relations, this identification permits us to define monic $\m$-modules and monomorphism category ${\rm mon}(B, \mathscr C) = {\rm mon}(Q, I, \mathscr C)$
([LZ1, LZ2], [ZX]), for any additive full subcategory $\mathscr C$ of $A$-mod.
This definition is combinatorial and constructive, and it admits a homological interpretation. In general, there is no longer the corresponding combinatorial definition of a monic module, but
this homological interpretation still makes sense, and it is taken as the definition of the monomorphism category ${\rm mon}(B, \mathscr C)$ by Hu, Luo, Xiong and Zhou in [HLXZ]. See Subsection 2.1.

\vskip5pt

The study of the monomorphism categories can be traced to G. Birkhoff [Bir].
When $B$ is the path algebra of quiver $A_n$ with linear orientation, i.e., $B = T_n(k) = \left(\begin{smallmatrix} k&k&\cdots & k \\ 0&k&\cdots &k \\  \vdots &\vdots & \vdots &\vdots \\0& 0 & \cdots & k\end{smallmatrix}\right)$, then $\m = A\otimes_kB = T_n(A)$ and the monomorphism category
${\rm mon}(B, A) = {\rm mon}(B, A\mbox{-}{\rm mod})$ is exactly {\it the submodule category} $\mathcal S_n(A)$ ([RS1-RS3]), or called {\it the filtered chain category} ([S1-S5]).
They have Auslander-Reiten sequences ([RS2]) and the RSS equivalence ([ZX]). Simson ([S2]-[S5]) has studied their representation type.
By  Kussin, Lenzing and Meltzer [KLM1, KLM2] and Chen [Chen1], they are related to the singularity theory.

\subsection{} An important application of monomorphism categories is that they can describe Gorenstein-projective modules as
${\rm gp}(\m) = {\rm mon}(B, \ {\rm gp}(A))$ (cf. Theorem \ref{gp} below). Thus,
Gorenstein-projective $\m$-modules are monic $\m$-modules over Gorenstein-projective $A$-modules.
If $B$ is given by a finite acyclic quiver and monomial relations, by the combinatorial definition of monic modules,
this gives in practice a reductive construction of Gorenstein-projective $\m$-modules.

\vskip5pt

{\bf Question 1.} \ Is there the similar relation between semi-Gorenstein-projective $\m$-modules and $A$-modules?

\begin{thm}\label{homodescriptionofsmon1}  Let $A$ and $B$ be finite-dimensional $k$-algebras with ${\rm gl.dim} B< \infty$, and $\m = A\otimes_kB$.
Then $${\rm mon}(B, \ ^\bot A) = {\rm mon}(B, A)\cap \ ^\bot \m.$$
Moreover, if \ ${\rm inj.dim} \ A_A < \infty$, then \ $^\bot \m =  {\rm mon}(B, \ ^\bot A)$.
\end{thm}

Theorem \ref{homodescriptionofsmon1} will be proved in Section 3, as a special case of Theorem \ref{homodescriptionofsmon}.

\subsection{} An Artin algebra $A$ is {\it Gorenstein}, if ${\rm inj.dim} \ _AA < \infty$ and ${\rm inj.dim} A_A < \infty.$ An Artin algebra $A$ is {\it left weakly Gorenstein} ([M1], [RZ2]), if any left semi-Gorenstein-projective $A$-module is Gorenstein-projective, i.e.,
$^\bot A = {\rm gp}(A).$ It is open whether a left weakly Gorenstein algebra is right  weakly Gorenstein ([M1, $\S$ 5], [RZ2, 9.3]).
However, if no confusions caused, we will omit the word ``left".

\vskip5pt

By Enochs and Jenda [EJ2, 11.5.3], Gorenstein algebras are weakly Gorenstein. By Yoshino [Y1, Theorem 5.5] and Beligiananis [Bel2, Corollary 5.11], if \ $^\perp A$ is of finite type, then $A$ is weakly Gorenstein. By
Marczinzik [M1, Theorem 3.5(3)], torsinless finite algebras are weakly Gorenstein.
For more information on weakly Gorenstein algebras we refer to [Bel1, Bel2], [M1], and [RZ2, 1.2 - 1.4, 3.6].

\vskip5pt

{\bf Question 2.} \ (i) \ Let $A$ and $B$ be Artin algebras, $M$ an $A$-$B$-bimodule such that $\m = \left(\begin{smallmatrix}A&M\\0&B\end{smallmatrix}\right)$ is an Artin algebra. When  $\m$ is weakly Gorenstein?

\vskip5pt

(ii) \  Let $A$ and $B$ be finite-dimensional $k$-algebras with ${\rm gl.dim} B< \infty$. When the tensor product $A\otimes_kB$ is weakly Gorenstein?

\vskip5pt

It turns out that, in the both cases, monic modules will play a crucial role.

\begin{thm}\label{wgandmon} \ Let $A$ and $B$ be Artin algebras, $M$ an $A$-$B$-bimodule with ${\rm proj.dim}_AM < \infty$, and $\m=\left(\begin{smallmatrix}A&M\\0&B\end{smallmatrix}\right)$.

\vskip5pt

$(1)$ \ If ${\rm proj.dim}M_B < \infty$ and ${\rm D}(M_B)\in \ (^\perp B)^{\perp}$,
then $\m$ is weakly Gorenstein if and only if each semi-Gorenstein-projective $\m$-module is monic respect to bimodule  $M$, and $A$ and $B$ are weakly Gorenstein.

\vskip5pt

$(2)$ \ If ${\rm proj.dim}M_B < \infty$ and $B$ is a Gorenstein algebra, then $\m$ is weakly Gorenstein if and only if each semi-Gorenstein-projective $\m$-module is monic respect to bimodule $M$ and $A$ is weakly Gorenstein.

\vskip5pt

$(3)$ \ If $_AM$ is torsionless and $M_B$ is projective, then $\m$ is weakly Gorenstein if and only if $A$ and $B$ are weakly Gorenstein.
In particular, if $_AM$ and $M_B$ are projective, then  $\m$ is weakly Gorenstein if and only if $A$ and $B$ are weakly Gorenstein.

\end{thm}

Theorem \ref{wgandmon} is the combination of Propositions \ref{wgandmon1},  \ref{cotorsionpair} and \ref{matrixwg1}.

\begin{thm}\label{tensorwg}  \ Let $A$ and $B$ be finite-dimensional $k$-algebras.

\vskip5pt

$(1)$ \ Assume that ${\rm gl.dim} B< \infty$ and $\m = A\otimes_kB$.
If $\m$ is weakly Gorenstein, then so is $A$. Conversely, if $A$ is weakly Gorenstein, then a semi-Gorenstein-projective $\m$-module $M$ is Gorenstein-projective if and only if $M$ is monic.

\vskip5pt

Thus, if $A$ is weakly Gorenstein,  then $\m$ is weakly Gorenstein if and only if each semi-Gorenstein-projective $\m$-module is monic,
or equivalently, \ $^\bot \m = {\rm mon}(B, \ ^\perp A)$.

\vskip5pt

$(2)$ \  Let $Q$ be a finite acyclic quiver. Then  $A\otimes_kkQ$ is weakly Gorenstein if and only if $A$ is weakly Gorenstein.
\end{thm}

The assumption ${\rm gl.dim} B< \infty$ holds automatically if $B$ is given by a bound acyclic quiver.
Theorems \ref {tensorwg} is the combination of Propositions \ref{wgandmon2} and \ref{quiverwg1}.

\subsection{} \ An Artin algebra $A$ will be called {\it left semi-Groenstein-projective-free}, or in short, lsgp-free, provided that each left semi-Gorenstein-projective $A$-module is a projective module, i.e.,  \ $^\perp A = \add(A)$. We do not know whether a lsgp-free algebra is right semi-Groenstein-projective-free.

\vskip5pt

Recall that $A$ is {\it left CM-free} ([Chen 2]) if \ ${\rm gp}(A) = \add(A)$. Thus, $A$ is lsgp-free if and only if
$A$ is left CM-free and left weakly Gorenstein. It is open whether a left CM-free algebra is left weakly Gorenstein (or equivalently, lsgp-free). See [RZ2, 9.2].
Many algebras are lsgp-free. For example, this is the case if ${\rm gl.dim.} A <\infty.$ There are also non Gorenstein algebras $A$ (thus, ${\rm gl.dim.} A = \infty$) which are lsgp-free.

\begin{thm}\label{lsgpfree}

\vskip5pt

$(1)$ \  Assume that \ $^\perp B = \add(B)$ with ${\rm proj.dim} M_B < \infty$ and that $_AM$ is torsionless with \ ${\rm proj.dim} \ _AM < \infty$.
Then $\m=\left(\begin{smallmatrix}A&M\\0&B\end{smallmatrix}\right)$ is weakly Gorenstein if and only if $A$ is weakly Gorenstein.

\vskip5pt

Moreover, if in addition $_AM$ is projective, then \ $$^\perp \m = {\rm gp}(\m) = \{\binom{M\otimes_BP}{P}_{{\rm Id}_{M\otimes_BP}} \oplus \binom{G}{0} \ | \ P\in \add(B), \ G\in \ ^\perp A = {\rm gp}(A)\}$$ and \ $^\perp \m = \add(\m)$ if and only if \ $^\perp A = \add(A)$.

\vskip5pt

$(2)$ \  Let $I$ be an admissible ideal of $kQ$, and $\m = A\otimes_kkQ/I$.
Then \ $^\perp \m = \add(\m)$ if and only if \ $^\perp A = \add(A)$.
\end{thm}

Theorems \ref {lsgpfree} is the combination of Propositions \ref{matrixwg2} and  \ref{quiverwg2}.

\subsection{} If $\m = A\otimes_kB$ with ${\rm gl.dim.} B<\infty$,
or if $\m = \left(\begin{smallmatrix}A & M \\ 0 & B\end{smallmatrix}\right)$, then Gorenstein-projective  $\m$-modules are always monic (cf. Theorems \ref{gp} and \ref{homodescriptionofsmon1}).
Is this true for semi-Gorenstein-projective $\m$-modules?  One may ask a stronger question:

\vskip5pt

{\bf Question 3.} \ In the both cases, whether there exist double semi-Gorenstein-projective $\m$-modules which are not monic?

\vskip5pt

The positive answer will in particular gives
double semi-Gorenstein-projective modules which are not Gorenstein-projective.
As mentioned in Subsection 1.2, this is highly nontrivial.

\vskip5pt

To answer {\bf Question 3}, we consider $\m = A\otimes_k k(\circ \longrightarrow \circ)=\left(\begin{smallmatrix}A & A \\ 0 & A\end{smallmatrix}\right) = T_2(A)$.
Any left $\m$-module $M$ can be identified with a triple $\left(\begin{smallmatrix}X\\ Y\end{smallmatrix}\right)_\varphi$, where $\varphi: Y \longrightarrow X$ is a left $A$-map.
Thus, one has the exact sequence of left $A$-modules
$Y \stackrel \varphi \longrightarrow X \stackrel \pi \longrightarrow {\rm Coker}\varphi \longrightarrow 0$,
and there is a unique $A$-map $\beta: {\rm Coker}\pi^* \longrightarrow Y^*$,  such that the diagram with exact rows
$$\xymatrix{0\ar[r] & ({\rm Coker}\varphi)^*\ar[r]^-{\pi^*}\ar@{=}[d] &  X^{*} \ar[r]^-{p}\ar@{=}[d] & {\rm Coker}\pi^*\ar@{-->}[d]^{\beta} \ar[r] & 0  \\
0 \ar[r] & ({\rm Coker}\varphi)^*\ar[r]^-{\pi^*} & X^* \ar[r]^-{\varphi^*} & Y^*.}$$
commutes, where $p$ is the canonical $A$-epimorphism. So one has the left $A$-map $\beta^*: Y^{**}\longrightarrow ({\rm Coker}\pi^*)^*$,  and the composition
$\beta^*\phi_Y: Y\longrightarrow ({\rm Coker}\pi^*)^*$, where $\phi_Y: Y\longrightarrow Y^{**}$ is the canonical map.

\begin{thm}\label{resultsforT2A} \ Let  $A$ be an Artin algebra,  $\m = T_2(A) = \left(\begin{smallmatrix}A&A\\0&A\end{smallmatrix}\right),$ and  $\left(\begin{smallmatrix}X\\ Y\end{smallmatrix}\right)_\varphi$
a left $\m$-module. Then

\vskip10pt

$(1)$ \ There is a left $\m$-module isomorphism \ $\left(\begin{smallmatrix}X\\ Y\end{smallmatrix}\right)_\varphi^{**}\cong  \left(\begin{smallmatrix}X^{**} \\ ({\rm Coker} \pi^*)^*\end{smallmatrix}\right)_{p^*}$,
where $p^{*}: ({\rm Coker} \pi^*)^*\longrightarrow X^{**}$ is the $A$-monomorphism induced by $p: X^*\longrightarrow ({\rm Coker}\varphi)^*$.

\vskip5pt

Taking this isomorphism as identity, then the canonical $\m$-map $\phi_{\binom{X}{Y}_\varphi}: \binom{X}{Y}_\varphi\longrightarrow \left(\begin{smallmatrix}X^{**} \\ ({\rm Coker} \pi^*)^*\end{smallmatrix}\right)_{p^*}$ is given by
$\phi_{\binom{X}{Y}_\varphi} = \binom{\phi_X}{\beta^*\phi_Y}.$

\vskip10pt

$(2)$ \ \ $\left(\begin{smallmatrix}X\\ Y\end{smallmatrix}\right)_\varphi$ is torsionless and double semi-Gorenstein-projective
if and only if  $\left(\begin{smallmatrix}X\\ Y\end{smallmatrix}\right)_\varphi$ is monic,
\ $X, \ Y$,  and  ${\rm Coker}\varphi$ are double semi-Gorenstein-projective, and $X$ and $Y$ are torsionless.

\vskip10pt

$(3)$ \ \ $\left(\begin{smallmatrix}X\\ Y\end{smallmatrix}\right)_\varphi$ is double semi-Gorenstein-projective
with epimorphism $\phi_{\left(\begin{smallmatrix}X\\ Y\end{smallmatrix}\right)_\varphi}$ if and only if \ $\varphi^*: X^*\longrightarrow Y^*$ is an epimorphism,
$X$ and $Y$ are double semi-Gorenstein-projective,  $({\rm Coker}\varphi)^*$ is semi-Gorenstein-projective, and \ $\phi_X$ and $\phi_Y$ are epimorphisms.
\end{thm}

Theorem \ref{resultsforT2A}(1) is a summary of Lemma \ref{doubledual} and Proposition \ref{thecanonicalmap};
and Theorem \ref{resultsforT2A}(2) and (3) will be clear after Proposition \ref{sgpsandref}.

\vskip5pt

As remarked in [RZ4, 3.1], up to now, all the known examples have the following property:

\vskip5pt

Double semi-Gorenstein-projective modules $M$ such that
$\phi_M$ is a monomorphism (an epimorphism, respectively) are Gorenstein-projective.

\vskip5pt

The following result shows that this property is preserved under the $T_2$-extensions.

\begin{thm} \label{sgpsandtorsionless} \ Let  $A$ be an Artin algebra and $\m = T_2(A) = \left(\begin{smallmatrix}A&A\\0&A\end{smallmatrix}\right).$
Then

\vskip5pt

$(\rm 1)$ \ \ Any torsionless and double semi-Gorenstein-projective $A$-module is Gorenstein-projective if and only if
any torsionless and double semi-Gorenstein-projective $\m$-module is Gorenstein-projective.

\vskip5pt

$(\rm 2)$ \ \ Any double semi-Gorenstein-projective $A$-module $L$ with $\phi_L$ an epimorphism is Gorenstein-projective
if and only if any double semi-Gorenstein-projective $\m$-module $M$ with $\phi_M$ an epimorphism
is Gorenstein-projective.
\end{thm}

Theorem \ref{sgpsandtorsionless} will be proved in Subsection 5.9.

\subsection{} The following result positively answers {\bf Question 3}, and
gives a construction of double semi-Gorenstein-projective $T_2(A)$-modules which are not monic.

\begin{thm} \label{doublesgpnotmonic} \ Suppose that $Y$ is a double semi-Gorenstein-projective $A$-module which is not torsionless.
Let $\varphi: Y\longrightarrow P$ be a left ${\rm add}(A)$-approximation of $Y$.  Then $\binom{P}{Y}_\varphi$ is a double semi-Gorenstein-projective
$T_2(A)$-module which is not monic. In particular, $\binom{P}{Y}_\varphi$ is not torsionless.
\end{thm}

Theorem \ref{doublesgpnotmonic} will be proved in Subsection 6.1.
Using the algebra $A$ in [RZ2] and the $A$-modules $M(1, -q, c)$ in
[RZ3], by Theorem \ref{doublesgpnotmonic}, we obtain a class of double semi-Gorenstein-projective $T_2(A)$-modules with parameter $c$ as
$$X(c): = \left(\begin{smallmatrix}_AA\\ M(1, -q, c)\end{smallmatrix}\right)_{f_1}$$
such that $X(c)$ not monic, and hence not torsionless; moreover, all the canonical  maps $\phi_{X(c)}: X(c)\longrightarrow X(c)^{**}$ are neither monomorphisms nor epimorphisms,
and $X(c)^{**}$ are not semi-Gorenstein-projective. See Proposition \ref{importantexamples}.

\vskip5pt

\section{\bf Preliminaries: Monic modules with relations to Gorenstein-projective modules}

\subsection{Monic modules over tensor algebras}

\begin{defn}\label{maindef1}  \ {\rm ([HLXZ, 3.1])} \ Let $A$ and $B$ be finite-dimensional $k$-algebras, and $\m = A\otimes_kB$.

\vskip5pt

{\rm (1)} \ \ A left $\m$-module $X$ is monic, if ${\rm Tor}^\m_i(A\otimes_k V, X) = 0$ for all $i\ge 1$ and for all right $B$-modules $V$.

\vskip5pt

Denote by ${\rm mon}(B, A)$ the full subcategory of $\m$-{\rm mod} consisting of monic modules, which is called the monomorphism category of $B$ over $A$.

\vskip5pt

{\rm (2)}  \ \ Let $\mathscr C$ be an additive full subcategory of $A$-{\rm mod}.
An object $X\in {\rm mon}(B, A)$ is a monic module over $\mathscr C$, if $(A\otimes_k V)\otimes_\m X\in \mathscr C$
for all right $B$-module $V$.
\vskip5pt

Denote by ${\rm mon}(B, \mathscr C)$ the full subcategory of ${\rm mon}(B, A)$ of monic modules  over $\mathscr C$, which is called the monomorphism category of $B$ over $\mathscr C$.
\end{defn}

\begin{lem}\label{desofsmon}  \  {\rm ([HLXZ, Lemma 3.2(7)]; [ZX, Theorem 2.6(1)])} \ One has \
\begin{align*}{\rm mon}(B, A) & = \{X\in\m\mbox{-}{\rm mod} \ | \ {\rm Tor}^\m_i(A\otimes_k {\rm D}(S), X) = 0, \ \forall \ i\ge 1, \ \forall \ \mbox{simple left} \ B\mbox{-module} \ S\} \\ & = \ ^\bot (\D(A_A)\otimes_k B).\end{align*}
\end{lem}

\begin{exm} \label{quiverexams} $(1)$  \ If $B$ is the path algebra of the quiver $A_n \ \ (n\ge 2)$ with linear orientation, then $B = T_n(k) = \left(\begin{smallmatrix} k&k&\cdots & k \\ 0&k&\cdots &k \\  \vdots &\vdots & \vdots &\vdots \\0& 0 & \cdots & k\end{smallmatrix}\right)$,  \ $\m =  A\otimes_k B = \left(\begin{smallmatrix} A&A&\cdots &A \\ 0&A&\cdots &A \\  \vdots &\vdots & \vdots &\vdots \\0& 0 & \cdots &  A\end{smallmatrix}\right) = T_n(A)$, and ${\rm mon}(B, A)$ turns out to be
$$\mathcal S_n(A) = \{\left(\begin{smallmatrix} X_1  \\  \vdots \\  \\ X_n\end{smallmatrix}\right)_{(\varphi_i)}\in T_n(A)\mbox{-}{\rm mod} \ | \ \varphi_i:
X_{i+1}\longrightarrow X_{i} \ \mbox{is a monomorphism}, \ \forall \ 1\le i\le n-1\}.$$
This submodule category has been studied in {\rm [A], [S1-S5], [RS1 - RS3], [Z1]}.

\vskip5pt

$(2)$  \ If $B$ is the path algebra $kQ$, where $Q = (Q_0, Q_1, s, e)$ is a finite acyclic quiver, then a monical $\m$-module has been
defined in {\rm [LZ1]} as a representation $(X_i, \ X_\alpha, \ i\in Q_0, \ \alpha\in Q_1)$ of $Q$ over $A$, such that  for each $i\in Q_0$ the $A$-map
$$(X_\alpha)_{\begin{smallmatrix} \alpha\in Q_1  \\  e(\alpha) = i\end{smallmatrix}}: \bigoplus_{\begin{smallmatrix} \alpha\in Q_1  \\  e(\alpha) = i\end{smallmatrix}}X_{s(\alpha)}\longrightarrow X_i$$
is a monomorphism.

\vskip5pt

For any additive full subcategory $\mathscr C$ of $A$-{\rm mod}, a monic $\m$-module over $\mathscr C$ has been defined in {\rm [ZX, 2.1]}, as a monic $\m$-module
$(X_i, \ X_\alpha, \ i\in Q_0, \ \alpha\in Q_1)$
satisfying $$X_i/{\rm Im}(X_\alpha)_{\begin{smallmatrix} \alpha\in Q_1  \\  e(\alpha)= i\end{smallmatrix}}\in \mathscr C, \ \ \forall \ i\in Q_0.$$

\vskip5pt

$(3)$ \ If $B = kQ/I$ with $I$ generated by monomial relations, then ${\rm mon}(B, \mathscr C)$ has
also been defined combinatorially. For details see {\rm [LZ2] and [ZX]}.

\vskip5pt

In all these monomorphism categories defined via quivers, ``monomorphisms" are visible,
and they also admit the homological description in {\rm Definition \ref {maindef1}} $(${\rm [Z1, Theorem 3.1], [LZ2, 2.1], [ZX, Theorem 2.6]}$)$.
\end{exm}

\begin{lem}\label{torsionlessmonic}  \ Let $\m = A\otimes_k kQ$, where $Q$ is a finite acyclic quiver.
Then torsionless $\m$-modules are monic.
\end{lem}
\noindent {\bf Proof.} \ Let $X = (X_i, \ X_\alpha)$ be a torsionless $\m$-module. Then $X$ is a submodule of a projective $\m$-module, which is of the form
$P\otimes_kL$, where $P$ is a projective left $A$-module, and $L = (L_i, \ L_\alpha)$ is a projective left $kQ$-module.  Thus there is a monomorphism
$(f_i)_{i\in Q_0}: (X_i, X_\alpha) \longrightarrow (P\otimes_k L(i), \ {\rm Id}_P\otimes_k L_\alpha)$ of $\Lambda$-modules. Hence, for each $i\in Q_0$, the diagram
of $A$-maps
$$\xymatrix{\bigoplus\limits_{\alpha\in Q_1, e(\alpha) = i} X_{s(\alpha)}
\ar@{^(->}[d]_-{\bigoplus\limits_{\alpha\in Q_1, e(\alpha)= i} f_{s(\alpha)}} \ar[rrr]^-{(X_\alpha)_{\alpha\in Q_1, e(\alpha) = i}} &&& X_i\ar@{^(->}[d]^-{f_i}
\\ \bigoplus\limits_{\alpha\in Q_1, e(\alpha) = i}P\otimes_k L_{s(\alpha)} \ar@{^(->}[rrr]^-{({\rm Id}_P\otimes_k L_\alpha)_{\alpha\in Q_1, e(\alpha) = i}}& & & P\otimes_k L_i}$$
commutes. Since both $\bigoplus\limits_{\alpha\in Q_1, e(\alpha)= i} f_{s(\alpha)}$ and $({\rm Id}_P\otimes_k L_\alpha)_{\alpha\in Q_1, e(\alpha) = i}$ are monomorphisms,
it follows that $(X_\alpha)_{\alpha\in Q_1, e(\alpha) = i}$ is a monomorphism, i.e., $X$ is a monic $\m$-module. $\s$

\vskip5pt

\begin{rem} \  {\rm Lemma \ref{torsionlessmonic}} is not true for $\m = A\otimes_k (kQ/I)$, even if $I$ is generated by monomial relations.
For example, if $A = k$, $Q = 3 \stackrel  \alpha \longrightarrow 2 \stackrel  \beta \longrightarrow 1$, and $I = \langle \beta\alpha\rangle$, then the simple module $S(2) = {\rm rad} P(3)$ is a torsionless
$(kQ/I)$-module, but it is not monic.
\end{rem}

\subsection{Gorenstein-projective modules over tensor algebras} The relationship between Gorenstein-projective modules over $\m = A\otimes_kB$
and monomorphism categories of $B$ over $A$ is

\begin{thm}\label{gp} {\rm ([HLXZ, Theorem 4.5])} \ Let $A$ and $B$ be finite-dimensional $k$-algebras with ${\rm gl.dim} B < \infty$, and $\m = A\otimes_kB$.
Then ${\rm gp}(\m) = {\rm mon}(B, \ {\rm gp}(A))$. In particular, a Gorenstein-projective $\m$-module is monic.
\end{thm}

Theorem \ref{gp} is proved for $\m = T_n(A) = A\otimes_kT_n(k)$ with $A$ Gorenstein in [Z1, Corollary 4.1(ii)];
it is proved for $B = kQ$ in [LZ1, Theorem 5.1], and for $B = kQ/I$ in [LZ2, Theorem 4.1], where $Q$ is any finite acyclic quiver, and $I$ is generated by monomial relations. In all these cases, since ${\rm mon}(B, \ {\rm gp}(A))$ are defined
via the combinatorics of quivers, Theorem \ref{gp} provides in practice an inductive construction of Gorenstein-projective modules.

\subsection{Monic modules respect to bimodules} Let $A$ and $B$ be Artin algebras, and $M$ an $A$-$B$-bimodule such that
$\Lambda=\left(\begin{smallmatrix}A&M\\0&B\end{smallmatrix}\right)$
is an Artin algebra. This is equivalent to say that $A$ and $B$ are Artin $R$-algebra, and $M$ is finitely generated over $R$ which acts centrally on $M$,
where $R$ is a commutative Artin ring ([ARS, Proposition 2.1, p.72]). Any left $\m$-module is identified with a triple  $\binom{X}{Y}_\varphi$,
where $X$ is a left $A$-module, \ $Y$ is a left $B$-module, and \ $\varphi: M\otimes_B Y \longrightarrow X$ is a left $A$-map.

\begin{defn}\label{maindef2} \ {\rm ([XZZ, 2.1])} \ Let $\Lambda=\left(\begin{smallmatrix}A&M\\0&B\end{smallmatrix}\right)$
be an Artin algebra. A $\m$-module $\binom{X}{Y}_\varphi$ is {\it monic respect to bimodule \ $_AM_B$},
provided that $\varphi: M\otimes_B Y \longrightarrow X$ is a monomorphism.

\vskip5pt Denote by $\mathcal M(A, M, B)$ the full subcategory of $\m$-mod of
monic $\m$-modules respect to bimodule \ $_AM_B$, which is called {\it the monomorphism category respect to  bimodule \ $_AM_B$}.
\end{defn}

\begin{exm} \label{bimoduleexams} \ The monomorphism category ${\rm mon}(B, A)$ and
the monomorphism category $\mathcal M(A, M, B)$ are in different setting.
Even if $\m=A\otimes_kB = \left(\begin{smallmatrix}A'&M\\0&B'\end{smallmatrix}\right)$, \ ${\rm mon}(B, A)\ne \mathcal M(A', M, B')$ in general.

\vskip5pt

For example, consider $T_n(A) =A\otimes_k T_n(k)$. A $T_n(A)$-module $X = \left(\begin{smallmatrix} X_1  \\  \vdots \\  \\ X_n\end{smallmatrix}\right)_{(\varphi_i)}$
is a monic $T_n(A)$-module if and only if $\varphi_i: X_{i+1}\longrightarrow X_{i}$ is a monomorphism for all \ $1\le i\le n-1.$  Thus ${\rm mon}(T_n(k), A) = \mathcal S_n(A)$.

\vskip5pt

On the other hand, $T_n(A) = \left(\begin{smallmatrix}T_{n-1}(A)&M_{n-1}\\0&A\end{smallmatrix}\right)$ for $n\ge 2$, where $M_{n-1} = \left(\begin{smallmatrix}A \\ \vdots \\ A\end{smallmatrix}\right)$ \ \   $($$n-1$ rows$)$
is a $T_{n-1}(A)$-$A$-bimodule, and $X$  is a monic $T_n(A)$-module respect to bimodule $M_{n-1}$ if and only if
$$\varphi_i\cdots\varphi_{n-1}: X_n \longrightarrow X_{i}$$ is a monomorphism for all \ $1\le i\le n-1.$

\vskip5pt

\vskip5pt

Thus, a $T_n(A)$-module $X$
is a monic $T_n(A)$-module if and only if
$\left(\begin{smallmatrix} X_1  \\  \vdots \\  \\ X_m\end{smallmatrix}\right)_{(\varphi_i)}$ is a monic $T_m(A)$-module respect to bimodule $M_{m-1}$ for all $2\le m\le n$, where $T_m(A) = \left(\begin{smallmatrix}T_{m-1}(A)&M_{m-1}\\0&A\end{smallmatrix}\right)$, and
$M_{m-1} = \left(\begin{smallmatrix}A \\ \vdots \\ A\end{smallmatrix}\right)$ \ \   $($$m-1$ rows$)$;
and a monic $T_n(A)$-module $X$ respect to $M_{n-1}$ is not necessarily a monic $T_n(A)$-module.
In some sense, $\mathcal M(T_{n-1}(A), M_{n-1}, A)$ can be view as the local version of ${\rm mon}(T_n(k), A)$.

For example, let $n\ge 3$. Consider $T_n(A)$-module $X = \left(\begin{smallmatrix} A  \\  \vdots \\ A \\ A\oplus A \\ X_n = A\end{smallmatrix}\right)_{(\varphi_i)}$, where
$$\varphi_{n-1} = \binom{{\rm Id}_A}{0}: A \longrightarrow A\oplus A, \ \ \ \
\varphi_{n-2} = ({\rm Id}_A, 0): A\oplus A \longrightarrow A$$ and $\varphi_i = {\rm Id}_A: X_{i+1}=A\longrightarrow A = X_{i}$ for all \ $1\le i\le n-3.$
Then $X\notin {\rm mon}(T_n(k), A)$, but $X\in \mathcal M(T_{n-1}(A), M_{n-1}, A)$.
\end{exm}

\begin{lem}\label{torsionlessmonic2}  \ Let $\Lambda=\left(\begin{smallmatrix}A&M\\0&B\end{smallmatrix}\right)$
be an Artin algebra, where $M_B$ is projective. Then torsionless $\m$-modules are monic respect to bimodule \ $_AM_B$.
\end{lem}
\noindent {\bf Proof.} \ Let $L = \binom{X}{Y}_\varphi$ be a torsionless $\m$-module. Then $L$ is a submodule of a projective $\m$-module, which is of the form
$\binom{P}{0}\oplus \binom{M\otimes_BQ}{Q}_{{\rm Id}_{M\otimes_BQ}}$, where $P$ is a projective left $A$-module, and $Q$ is a projective left $B$-module.
Thus, there is a monomorphism $\binom{\binom{f_1}{f_2}}{g}: \binom{X}{Y}_\varphi\longrightarrow \binom{P\oplus M\otimes_BQ}{Q}_{\binom{0}{1}}.$ Since $M_B$ is projective, ${\rm Id}_M\otimes_Bg$ is a monomorphism.
By the commutative diagram
$$\xymatrix{M\otimes_BY \ar@{^(->}[d]_-{{\rm Id}_M\otimes_Bg} \ar[r]^-{\varphi} & X\ar@{^(->}[d]^-{\binom{f_1}{f_2}}
\\ M\otimes_BQ\ar@{^(->}[r]^-{\binom{0}{1}}& P\oplus M\otimes_BQ}$$
$\varphi$ is a monomorphism, i.e., $L$ is a monic $\m$-module. $\s$

\vskip5pt

\begin{rem} \  {\rm Lemma \ref{torsionlessmonic2}} is not true if $M_B$ is not projective.
For example, let $\Lambda=\left(\begin{smallmatrix}k&M\\0&A\end{smallmatrix}\right)$, where $A$ is the path algebra $k(2\longrightarrow 1)$, \
$_kM_A = \D(Ae_1) = \Hom_k(Ae_1, k)$ is a $k$-$A$-bimodule. Since $M e_2 = 0$, $M\otimes_A Ae_2 = M e_2\otimes_A e_2 = 0.$
Thus $\binom{0}{Ae_2}$ is a left projective $\m$-module. Let $\sigma: Ae_1 \longrightarrow Ae_2$ be the embedding. Then
$\binom{0}{\sigma}: \ \binom{0}{Ae_1} \longrightarrow \binom{0}{Ae_2}$ \ is a $\m$-monomorphism, and hence $\binom{0}{Ae_1}$ is a torsionless $\m$-module.
But since $M\otimes_A Ae_1= M\otimes_A e_1Ae_1 = k\ne 0$,  $\binom{0}{Ae_1}$ is not monic respect to bimodule \ $_kM_A$.
\end{rem}

\subsection{(Semi-)Gorenstein-projective modules over triangular matrix algebras of bimodules} For an Artin algebra $B$, let ${\rm D}$ denote the duality of $B$ ([ARS, p.37]).

\vskip5pt

Let $A$ and $B$ be Artin algebras, $M$ an $A$-$B$-bimodule such that $\Lambda=\left(\begin{smallmatrix}A&M\\0&B\end{smallmatrix}\right)$ is an Artin algebra.
Under suitable conditions, semi-Gorenstein-projective $\m$-modules can be described as follows.

\vskip5pt

\begin{thm} \label{sgpovermatrixalg} \   {\rm ([XZ, Theorem 1.1])} \
Assume that ${\rm proj.dim}_AM < \infty$ and ${\rm D}(M_B)\in \ (^\perp (_BB))^{\perp}$. Then a $\m$-module $\left(\begin{smallmatrix}X\\Y\end{smallmatrix}\right)_{\varphi}\in \
^\perp \m$ if and only if \ $Y\in \ ^\perp B$, the left $A$-map \ $\varphi: M\otimes_BY\longrightarrow X$
induces isomorphisms $\Ext^i_A(X, A)\cong\Ext^i_A(M\otimes_BY, A)$ for all $i\geq1$, and $\varphi^*: X^*\longrightarrow (M\otimes_BY)^*$ is a right $A$-epimorphism.
\end{thm}

An  $A$-$B$-bimodule $M$ is {\it compatible} ([Z2, Definition 1.1]),
if the following two conditions hold:

\vskip5pt

If $Q^\bullet$ is an exact sequence of projective
$B$-modules, then $M\otimes_BQ^\bullet$ is exact; and

\vskip5pt

If $P^\bullet$ is a complete $A$-projective
resolution, then $\Hom_A(P^\bullet, M)$ is exact.

\vskip5pt

\begin{lem} \label{compatible} \  {\rm ([Z2, Proposition 1.3(1)])} \ Let $M$ be an
$A$-$B$-bimodule. If ${\rm proj.dim} _AM < \infty$ and ${\rm proj.dim} M_B < \infty$, then $M$ is compatible.
\end{lem}

Under the condition of compatible bimodule, Gorenstein-projective $\m$-modules can be described as follows. In particular,  again, Gorenstein-projective $\m$-modules are monic, but in the sense of respect to bimodule $_AM_B$ (compare Theorem \ref{gp}).

\begin{thm} \label{gpovermatrixalg} \ {\rm ([Z2, Theorem 1.4])} \  Assume that $M$ is a compatible $A$-$B$-bimodule. Then $\left(\begin{smallmatrix}X\\Y\end{smallmatrix}\right)_{\phi}\in \
{\rm gp}(\m)$ if and only if  $\varphi: M\otimes_BY \longrightarrow X$ is a monomorphism,  $\cok\varphi\in {\rm gp}(A)$, and
$Y\in \ {\rm gp} (B)$.

If ths is the case, $X\in  {\rm gp}(A)$ if and only if $M\otimes_BY\in {\rm gp}(A).$
\end{thm}

\begin{cor} \label{sgpT2A} \  Let $A$ be an Artin algebra, and $\m = T_2(A) =\left(\begin{smallmatrix}A&A\\0&A\end{smallmatrix}\right)$.
Then

\vskip8pt

$(1)$  \ \ $^\perp \m  = \{\left(\begin{smallmatrix}X\\Y\end{smallmatrix}\right)_{\varphi}\in \m\mbox{-}{\rm mod} \ | \ X\in \ ^\perp A, \ \ Y\in \ ^\perp A, \ \ \varphi^*: X^*\longrightarrow Y^* \ \mbox{is epic}\}.$

\vskip8pt

$(2)$ \ \ ${\rm gp}(\m)   = \{\left(\begin{smallmatrix}X\\Y\end{smallmatrix}\right)_{\varphi}\in \m\mbox{-}{\rm mod} \ | \ \varphi: Y \longrightarrow X \ \mbox{is monic}, \ \
\cok\varphi\in {\rm gp}(A), \  \ Y\in {\rm gp} (A)\}$

\vskip5pt\hskip50pt $= \{\left(\begin{smallmatrix}X\\Y\end{smallmatrix}\right)_{\varphi}\in \m\mbox{-}{\rm mod} \ | \ \varphi: Y \longrightarrow X \ \mbox{is monic}, \ \
\cok\varphi\in {\rm gp}(A), \  \ Y\in {\rm gp} (A), \ \ X\in {\rm gp}(A)\}.$
\end{cor}

\vskip5pt

\section{\bf Monomorphism categories over perpendicular categories}

Let $A$ and $B$ be finite-dimensional $k$-algebras, and $\m = A\otimes_kB$.
A relation between semi-Gorenstein-projective $\m$-modules and semi-Gorenstein-projective $A$-modules
is contained in the following general result.

\vskip5pt

\begin{thm}\label{homodescriptionofsmon}  Let $A$ and $B$ be finite-dimensional $k$-algebra with ${\rm gl.dim} B< \infty$, $T$ an $A$-module, and $\m = A\otimes_kB$.
Then \ $${\rm mon}(B, \ ^\bot T) = {\rm mon}(B, A)\cap \ ^\bot (T\otimes_kB).$$

\vskip5pt

Moreover, if there is an exact sequence of left $A$-modules $$0\longrightarrow T_m\longrightarrow
\cdots\longrightarrow T_0\longrightarrow \D(A_A)\longrightarrow 0$$ with each
$T_j\in {\rm add} (T)$, then \ ${\rm mon}(B, \ ^\bot T) = \ ^\bot (T\otimes_kB)$.

\vskip5pt

In particular, there holds \ ${\rm mon}(B, \ ^\bot A) = {\rm mon}(B, A)\cap \ ^\bot \m;$ and if \ ${\rm inj.dim} \ A_A < \infty$, then \ $^\bot \m =  {\rm mon}(B, \ ^\bot A)$.
\end{thm}
\noindent{\bf Proof.} \ Let $X\in {\rm mon}(B, A)$. Since by definition ${\rm mon}(B, \ ^\bot T) \subseteq {\rm mon}(B, A),$ it follows that, in order to prove
${\rm mon}(B, \ ^\bot T) = {\rm mon}(B, A)\cap \ ^\bot (T\otimes_kB),$  it suffices to prove that  $X\in \ ^\bot (T\otimes_kB)$  if and only if $X\in {\rm mon}(B, \ ^\bot T)$, i.e.,
$(A\otimes_k V)\otimes_\m X\in \ ^\bot T$ for all right $B$-modules $V$.

\vskip5pt

Take a $\m$-projective resolution $$P_\bullet: \ \ \cdots \longrightarrow
P_1 \longrightarrow P_0 \longrightarrow
X\longrightarrow 0.$$

{\bf Claim 1:} \ $X\in \ ^\bot (T\otimes_kB)$ if and only if the complex $$\Hom_\m(P_\bullet, \Hom_k(S', T))$$ is exact, for each right simple $B$-module $S'$.

\vskip5pt

Since ${\rm gl.dim} B< \infty$, it is clear that $^\bot (T\otimes_kB) = \bigcap\limits_{S}\ ^\bot (T\otimes_kS)$, where $S$ ranges over all the left simple $B$-modules.
To use the Tensor-$\Hom$ adjoint pair later, we write a left simple $B$-module $S$ as $\D(S')$, where $S'$ is a right simple $B$-module.
Thus, $^\bot (T\otimes_kB) = \bigcap\limits_{S'}\ ^\bot (T\otimes_k \D(S'))$, where $S'$ ranges over all the right simple $B$-modules.
Therefore,  $X\in \ ^\bot (T\otimes_kB)$ if and only if $$\Hom_\m(P_\bullet, T\otimes_k\D(S'))$$ is exact, for each right simple $B$-module $S'$.
Note that the canonical $k$-linear isomorphism
$$T\otimes_k\D(S') \longrightarrow \Hom_k(S', T), \ t\otimes f\mapsto ``s'\mapsto f(s')t", \ \forall \ t\in T, \ f\in \D(S'), \ s'\in S'$$
is a left $\m$-isomorphism. Thus,
$X\in \ ^\bot (T\otimes_kB)$ if and only if \ $\Hom_\m(P_\bullet, \Hom_k(S', T))$ is exact, for each right simple $B$-module $S'$.

\vskip5pt

{\bf Claim 2:} \ $(A\otimes_k V)\otimes_\m X\in \ ^\bot T$ for all right $B$-modules $V$ if and only if the complex
$$\Hom_A((A\otimes_k S')\otimes_\m P_\bullet, T)$$ is exact for each right simple $B$-module $S'$.

\vskip5pt

By assumption $X\in {\rm mon}(B, A)$, i.e.,
${\rm Tor}^\m_i(A\otimes_k V, X) = 0$ for all $i\ge 1$ and for all right $B$-modules $V$. It follows that the functor
$$(A\otimes_k -)\otimes_\m X: {\rm mod}B \longrightarrow A\mbox{-}{\rm mod}$$ is an exact functor. As a consequence,
$(A\otimes_k V)\otimes_\m X\in \ ^\bot T$ for all right $B$-modules $V$ if and only if
$(A\otimes_k S')\otimes_\m X\in \ ^\bot T$ for each right simple $B$-module $S'$, since $^\bot T$ is extension closed.
Since ${\rm Tor}^\m_i(A\otimes_k S', X) = 0$ for all $i\ge 1$, it follows that
$$(A\otimes_kS')\otimes_\m P_\bullet: \ \ \cdots \longrightarrow (A\otimes_kS')\otimes_\m P_1 \longrightarrow \cdots \longrightarrow (A\otimes_kS')\otimes_\m P_0 \longrightarrow (A\otimes_kS')\otimes_\m X\longrightarrow 0$$
is an exact sequence of left $A$-modules. Since each $P_i$ is a projective left $\m$-module,
each  $(A\otimes_kS')\otimes_\m P_i\in \add (A\otimes_kS').$
Thus each  $(A\otimes_kS')\otimes_\m P_i$ is projective as a left $A$-module, and hence $(A\otimes_kS')\otimes_\m P_\bullet$  is an $A$-projective resolution of left $A$-module $(A\otimes_kS')\otimes_\m X$, for each
right simple $B$-module $S'$.
Therefore, $(A\otimes_k S')\otimes_\m X\in \ ^\bot T$ for each right simple $B$-module $S'$ if and only if
$\Hom_A((A\otimes_k S')\otimes_\m P_\bullet, T)$ is exact for each right simple $B$-module $S'$.

\vskip5pt

{\bf Claim 3:} \ There is an isomorphism of complexes
$$\Hom_A((A\otimes_k S')\otimes_\m P_\bullet, T)\cong \Hom_\m(P_\bullet, \Hom_k(S', T))$$
for each right simple $B$-module $S'$.

\vskip5pt

Applying the Tensor-$\Hom$ adjoint pair $((A\otimes_k S')\otimes_\m -, \ \Hom_A(A\otimes_k S', -))$ between $\m$-mod and $A$-mod, one has the following isomorphism of complexes of $k$-spaces
$$\Hom_A((A\otimes_k S')\otimes_\m P_\bullet, T)\cong \Hom_\m(P_\bullet, \Hom_A(A\otimes_k S', T)).$$
Applying the adjoint pair $(A\otimes_k -, \ \Hom_A(A, -))$ between $k$-mod and $A$-mod, one has the isomorphisms of $k$-spaces
$$\Hom_A(A\otimes_k S', T)\cong \Hom_k(S', \Hom_A(A, T))\cong \Hom_k(S', T),$$
which is clearly also an isomorphism of left $\m$-modules. All together we get
an isomorphism of complexes
$$\Hom_A((A\otimes_k S')\otimes_\m P_\bullet, T)\cong \Hom_\m(P_\bullet, \Hom_k(S', T))$$
for each right simple $B$-module $S'$.

\vskip5pt

It follows from {\bf Claim 1}, {\bf Claim 2} and {\bf Claim 3} that  $X\in \ ^\bot (T\otimes_kB)$  if and only if $(A\otimes_k V)\otimes_\m X\in \ ^\bot T$ for all right $B$-module $V$.
This proves \ ${\rm mon}(B, \ ^\bot T) = {\rm mon}(B, A)\cap \ ^\bot (T\otimes_kB)$.

\vskip5pt

Finally, assume that there is an exact sequence $0\longrightarrow T_m\longrightarrow
\cdots\longrightarrow T_0\longrightarrow \D(A_A)\longrightarrow 0$ with each
$T_j\in {\rm add} (T)$. To show ${\rm mon}(B, \ ^\bot T)= \
^\bot (T\otimes_k B)$, it suffices to show
$^\bot (T\otimes_k B)\subseteq \ {\rm mon}(B, A)$.
By Lemma \ref{desofsmon}, ${\rm mon}(B, A) = \ ^\bot (\D(A_A)\otimes_k B)$. Thus,
it suffices to show
$^\bot (T\otimes_k B)\subseteq \ ^\bot (\D(A_A)\otimes_k B)$.
This
follows from the exact sequence \ $0\longrightarrow T_m\otimes_k B
\longrightarrow \cdots \longrightarrow T_0\otimes_k B\longrightarrow
\D(A_A)\otimes_k B \longrightarrow 0$ \  with each $T_j \otimes_k B\in
\mbox{add}(T\otimes_k B).$
This completes the proof.  $\s$

\vskip5pt

\section{\bf Weakly Gorenstein algebras: Proof of Theorems \ref{wgandmon}, \ref{tensorwg}, and \ref{lsgpfree}}

\subsection{When triangular matrix algebras of bimodules are weakly Gorenstein?} Let $A$ and $B$ be Artin algebras, $M$ an $A$-$B$-bimodule such that $\Lambda=\left(\begin{smallmatrix}A&M\\0&B\end{smallmatrix}\right)$ is an Artin algebra. We will give various conditions for $\m$ being a left weakly Gorenstein algebra, i.e., $\ ^\perp \m = {\rm gp}(\m)$.

\begin{prop}\label{wgandmon1} \ Assume that ${\rm proj.dim}_AM < \infty$,  ${\rm proj.dim}M_B < \infty$, and ${\rm D}(M_B)\in \ (^\perp (_BB))^{\perp}$.
Then $\m =\left(\begin{smallmatrix}A&M\\0&B\end{smallmatrix}\right)$ is weakly Gorenstein if and only if each semi-Gorenstein-projective $\m$-module is monic respect to bimodule \ $_AM_B$, and $A$ and $B$ are weakly Gorenstein.
\end{prop}

\noindent{\bf Proof.}  \ Since ${\rm proj.dim}_{A}M <\infty$ and ${\rm proj.dim}M_{B} <\infty$,
the $A$-$B$-bimodule $M$ is compatible (cf. Lemma \ref{compatible}). Thus,  under the assumptions, one can apply Theorems \ref{sgpovermatrixalg} and \ref{gpovermatrixalg}.

\vskip5pt

Assume that each semi-Gorenstein-projective $\m$-module is monic respect to bimodule \ $_AM_B$, and that $A$ and $B$ are weakly Gorenstein. Let $\binom{X}{Y}_\varphi\in \ ^\perp \m.$
We need to prove $\binom{X}{Y}_\varphi\in \ {\rm gp}(\m).$ By the assumption, $\varphi: M\otimes_BY\longrightarrow X$ is a monomorphism; together with Theorem \ref{sgpovermatrixalg}, one gets the conclusions:

\vskip5pt

$\bullet$ \ $\varphi: M\otimes_BY\longrightarrow X$ is a monomorphism;

$\bullet$  \ $Y\in \ ^\perp B,$ and hence $Y\in {\rm gp}(B)$ \ (since by assumption $B$ is weakly Gorenstein);

$\bullet$ \ $\varphi$ induces isomorphisms $\Ext^i_A(X, A)\cong\Ext^i_A(M\otimes_BY, A)$  for all $i\geq1$;

$\bullet$ \ $\varphi^*: X^*\longrightarrow (M\otimes_BY)^*$ is a right $A$-epimorphism.

\vskip5pt

Applying $\Hom_A(-, A)$ to the exact sequence $0\longrightarrow M\otimes_B Y \stackrel \varphi \longrightarrow X \longrightarrow \cok \varphi \longrightarrow 0$,
since $\varphi^*: X^*\longrightarrow (M\otimes_BY)^*$ is an epimorphism and  $\varphi$ induces isomorphisms $\Ext^i_A(X, A)\cong\Ext^i_A(M\otimes_BY, A)$ for all $i\geq1$,
it follows that $\cok \varphi\in \ ^\bot A$.  Hence $\cok \varphi\in {\rm gp}(A)$, since by assumption $A$ is weakly Gorenstein.
Thus, we get the following:

\vskip5pt

$\bullet$ \  $\varphi: M\otimes_B Y \longrightarrow X$ is a monomorphism;

$\bullet$ \ $\cok \varphi\in \ {\rm gp}(A)$; and

$\bullet$ \ $Y\in {\rm gp}(B)$.

\vskip5pt

\noindent Applying Theorem \ref{gpovermatrixalg}, one gets $\binom{X}{Y}_\varphi \in {\rm gp}(\m)$. This proves the ``if" part.

\vskip10pt

Conversely, assume that $\m$ is weakly Gorenstein. Thus, any semi-Gorenstein-projective $\m$-module is Gorenstein-projective, and hence it is monic respect to bimodule \ $_AM_B$,
by Theorem \ref{gpovermatrixalg}. It remains to prove that $A$ and $B$ are weakly Gorenstein.
Let $X\in \ ^\perp A$. Applying  Theorem \ref{sgpovermatrixalg} one knows $\binom{X}{0}\in \ ^\perp \m$, thus
$\binom{X}{0}\in {\rm gp}(\m)$ by the assumption, and then by Theorem \ref{gpovermatrixalg} one has $X\in {\rm gp}(A).$ This proves that $A$ is weakly Gorenstein.

Similarly, let $Y\in \ ^\perp B$. By  Theorem \ref{sgpovermatrixalg} one knows $\binom{M\otimes_BY}{Y}_{{\rm Id}_{M\otimes_BY}}\in \ ^\perp \m$, and hence
$\binom{M\otimes_BY}{Y}_{{\rm Id}_{M\otimes_BY}}\in {\rm gp}(\m).$ Then by Theorem \ref{gpovermatrixalg}, $Y\in {\rm gp}(B).$ This proves that $B$ is weakly Gorenstein.
$\s$

\vskip10pt

Taking $B$ to be a Gorenstein algebra in Proposition \ref{wgandmon1}, we get

\vskip5pt

\begin{prop}\label{cotorsionpair} \ Assume that ${\rm proj.dim}_AM < \infty$,  ${\rm proj.dim}M_B < \infty$, and that $B$ is a Gorenstein algebra.
Then $\m =\left(\begin{smallmatrix}A&M\\0&B\end{smallmatrix}\right)$ is weakly Gorenstein if and only if each semi-Gorenstein-projective $\m$-module is monic respect to bimodule \ $_AM_B$ and $A$ is weakly Gorenstein.
\end{prop}

\noindent{\bf Proof.}  \ Since by assumption $B$ is a Gorenstein algebra, it follows that $\ ^\perp (_BB) = {\rm gp}(B)$. Recall that
for  a Gorenstein algebra $B$, $({\rm gp}(B), {\rm p}(B)^{< \infty})$ is a cotorsion pair (see e.g., [H], [EJ2], [BR]), where ${\rm p}(B)^{< \infty}$ is the full subcategory of $B$-mod consisting of modules of finite projective dimension. Thus $$\ (^\perp (_BB))^{\perp} = {\rm gp}(B)^{\perp} = {\rm p}(B)^{< \infty}.$$
Since by assumption ${\rm proj.dim}M_B < \infty$, It follows that  ${\rm inj.dim} {\rm D}(M_B) < \infty$. Since $B$ is Gorenstein,  it follows that ${\rm proj.dim} {\rm D}(M_B) < \infty$, i.e.,
${\rm D}(M_B)\in {\rm p}(B)^{< \infty} = \ (^\perp (_BB))^{\perp}$. Thus, the assertion follows from Proposition \ref{wgandmon1}. $\s$

\vskip5pt

Taking $B$ to be a field $k$ in Proposition \ref{cotorsionpair}, we get

\vskip5pt

\begin{cor}\label{corctp} \ Let $A$ be a finite-dimensional $k$-algebra.

\vskip5pt

$(1)$ \ Let $M$ be a finite-dimensional $A$-module. Assume that ${\rm proj.dim}_AM < \infty$.
Then $\m =\left(\begin{smallmatrix}A& M\\0& k\end{smallmatrix}\right)$ is weakly Gorenstein if and only if each semi-Gorenstein-projective $\m$-module is monic respect to bimodule \ $_AM_k$ and $A$ is weakly Gorenstein.

\vskip5pt

$(2)$ \ Let $P$ be a finite-dimensional projective left $A$-module, and $\m =\left(\begin{smallmatrix}A& P\\0& k\end{smallmatrix}\right)$.
Then
$$^\perp \m = \{\binom{G}{0} \oplus \binom{P\otimes_kV}{V}_{{\rm Id}_{P\otimes_kV}} \ | \ G\in \ ^\perp A, \ V\in k\mbox{-}{\rm mod}\}$$
$${\rm gp}(\m) = \{\binom{G}{0} \oplus \binom{P\otimes_kV}{V}_{{\rm Id}_{P\otimes_kV}} \ | \ G\in \ {\rm gp}(A), \ V\in k\mbox{-}{\rm mod}\}$$

\vskip5pt

\noindent and $\m$ is weakly Gorenstein if and only if $A$ is weakly Gorenstein.
\end{cor}

\noindent{\bf Proof.}  \ $(2)$  \ Let $\binom{X}{V}_\varphi\in \ ^\perp \m.$ By Theorem \ref{sgpovermatrixalg}, one has

$\bullet$ \ $\varphi: P^{\oplus {\rm dim}V} \longrightarrow X$
induces isomorphisms $\Ext^i_A(X, A)\cong\Ext^i_A(P^{\oplus {\rm dim}V}, A) = 0, \ \forall \ i\geq1$; and

$\bullet$ \ $\varphi^*: X^*\longrightarrow (P^*)^{\oplus {\rm dim}V}$ is a right $A$-epimorphism.

\vskip5pt
\noindent Thus $X\in \ ^\perp A$, and $\varphi^*$ is a splitting epimorphism. Hence
$\varphi^{**}$ is a splitting monomorphism. By the commutative diagram
$$\xymatrix{
P^{\oplus {\rm dim}V}\ar[d]_{\cong} \ar[r]^-{\varphi} & X \ar[d]^-{\phi_{X}}\\
(P^{**})^{\oplus {\rm dim}V}\ar[r]^-{\varphi^{**}} & X^{**}}$$
one sees that $\varphi$ is a splitting monomorphism. Thus $X = G\oplus P^{\oplus {\rm dim}V}$ where $G$ is semi-Gorenstein-projective, and hence
$\binom{X}{V}_\varphi \cong \binom{G}{0}\oplus \binom{P\otimes_k V}{V}_{\rm Id}$. This proves $^\perp \m = \{\binom{G}{0} \oplus \binom{P\otimes_kV}{V}_{{\rm Id}_{P\otimes_kV}} \ | \ G\in \ ^\perp A, \ V\in k\mbox{-}{\rm mod}\}$.
Since $\binom{P\otimes_k V}{V}_{\rm Id}$ is a projective $\m$-module and $\binom{G}{0}$  is a Gorenstein-projective $\m$-module if and only if $G$ a Gorenstein-projective $A$-module, it follows that
${\rm gp}(\m) = \{\binom{G}{0} \oplus \binom{P\otimes_kV}{V}_{{\rm Id}_{P\otimes_kV}} \ | \ G\in \ {\rm gp}(A), \ V\in k\mbox{-}{\rm mod}\}$.
Therefore $\m$ is weakly Gorenstein if and only if $A$ is weakly Gorenstein. $\s$

\vskip10pt

\begin{prop}\label{matrixwg1} \ Assume that $_AM$ is torsionless with \ ${\rm proj.dim} \ _AM < \infty$ and that $M_B$ is projective.
Then $\m=\left(\begin{smallmatrix}A&M\\0&B\end{smallmatrix}\right)$ is weakly Gorenstein if and only if $A$ and $B$ are weakly Gorenstein.

\vskip5pt

In particular, if $_AM$ and $M_B$ are projective, then  $\m$ is weakly Gorenstein if and only if $A$ and $B$ are weakly Gorenstein. \end{prop}

\noindent{\bf Proof.}  \ Since $M_B$ is projective, $\D(M_B)$ is an injective $B$-module, and hence $\D(M_B)\in (^\perp B)^\perp$.
Thus, the assumption that ${\rm proj.dim.} _AM <\infty$ and $M_B$ is projective guarantee that the conditions of Proposition \ref{wgandmon1} are satisfied.
By Proposition \ref{wgandmon1}, it suffices to prove that if $A$ and $B$ are weakly Gorenstein, then any semi-Gorenstein-projective $\m$-module $\binom{X}{Y}_\varphi$ is monic respect to bimodule \ $_AM_B$.

\vskip5pt

Applying Theorem \ref{sgpovermatrixalg} to $\binom{X}{Y}_\varphi\in \ ^\perp \m,$ one gets the following conclusions:

\vskip5pt

$\bullet$  \ $Y\in \ ^\perp B,$ and hence $Y\in {\rm gp}(B)$ \ (since by assumption $B$ is weakly Gorenstein);

$\bullet$ \ $\varphi: M\otimes_BY\longrightarrow X$
induces isomorphisms $\Ext^i_A(X, A)\cong\Ext^i_A(M\otimes_BY, A), \ \forall \ i\geq1$; and

$\bullet$ \ $\varphi^*: X^*\longrightarrow (M\otimes_BY)^*$ is a right $A$-epimorphism.

\vskip5pt

Since $Y\in {\rm gp}(B),$ \  $_BY$ is a submodule of some projective $B$-module $_BP$. Since $M_B$ is projective, it follows that
\ $_A(M\otimes_BY)$ is a submodule of \ $_A(M\otimes_BP)$. Since $_BP$ is projective,
\ $_A(M\otimes_BP)\in \add(_AM)$. Since by assumption $_AM$ is torsionless,
it follows that $M\otimes_BP$ is a torsionless left $A$-module, and hence $M\otimes_BY$ is a torsionless left $A$-module. Thus, the canonical map
$\phi_{M\otimes_BY}: M\otimes_BY \longrightarrow (M\otimes_BY)^{**}$ is a monomorphism.

\vskip5pt

Since $\varphi^*: \ X^* \longrightarrow (M\otimes_B Y)^*$ is an epimorphism, it follows that
$\varphi^{**}: \ (M\otimes_B Y)^{**} \longrightarrow X^{**}$ is a monomorphism. From the commutative diagram with monomorphism $\phi_{M\otimes_BY}$
$$\xymatrix{
M\otimes_B Y\ar@{^{(}->}[d]_{\phi_{M\otimes_B Y}} \ar[r]^-{\varphi} & X \ar[d]^-{\phi_{X}} \\
(M\otimes_B Y)^{**} \ar@{^{(}->}[r]^-{\varphi^{**}} & X^{**}}$$
one sees that $\varphi: M\otimes_B Y \longrightarrow X$ is a monomorphism, i.e., $\binom{X}{Y}_\varphi$ is monic respect to bimodule \ $_AM_B$. This completes the proof.
$\s$

\begin{rem} \ The ``only if" part in {\rm Proposition \ref{matrixwg1}} does not need the  condition that $_AM$ is torsionless.
\end{rem}

\subsection{When tensor algebras are weakly Gorenstein?} \ Let $A$ and $B$ be finite-dimensional $k$-algebras with ${\rm gl.dim} B< \infty$, and $\m = A\otimes_kB$. We first look at some properties of a map ${\rm mon}(B, -)$.

\begin{lem}\label{monmap}  \ Let $A$ and $B$ be finite-dimensional $k$-algebra, and $\m = A\otimes_kB$.

\vskip5pt

${\rm (i)}$ \ Let $\mathscr C$  be an additive full subcategory of $A$-mod closed under direct summands, and $M\in A$-mod.
Then $M\otimes_kB\in {\rm mon}(B, \mathscr C)$ if and only if $M\in \mathscr C.$

\vskip5pt

${\rm (ii)}$  \ If $M$ is a semi-Gorenstein-projective $A$-module  which is not Gorenstein-projective, then
$M\otimes_kB$ is a semi-Gorenstein-projective $\m$-module which is not Gorenstein-projective.

\vskip5pt

${\rm (iii)}$  \ Let \ $\Omega$ $($respectively, $\Gamma$$)$ be the set of additive full subcategories of $A$-mod $($respectively, $\m$-mod$)$ closed under direct summands. Then the map $${\rm mon}(B, -): \Omega \longrightarrow \Gamma, \ \ \ \mathscr C\mapsto {\rm mon}(B, \mathscr C)$$ is an injective map.

\end{lem}
\noindent{\bf Proof.} \ ${\rm (i)}$ \ Assume that $M\in \mathscr C.$ For any right $B$-module $V$, taking a $B$-projective resolution
$$P_\bullet: \ \cdots \longrightarrow P_1\longrightarrow P_0 \longrightarrow V \longrightarrow 0$$
of $V$, one has a projective resolution $A\otimes_kP_\bullet$ of right $\m$-module $A\otimes_kV$.  By the isomorphisms
$$(A\otimes_k V)\otimes_\m(M\otimes_k B)\cong (A\otimes_A M)\otimes_k(V\otimes_B B)\cong M\otimes_kV$$ one sees that there is an isomorphism of complexes
$$(A\otimes_kP_\bullet)\otimes_\m (M\otimes_k B)\cong M\otimes_k P_\bullet$$ and hence $${\rm Tor}^\m_i(A\otimes_k V, M\otimes_k B) \cong {\rm Tor}^k_i(M, V) = 0.$$
This shows $M\otimes_k B\in {\rm mon}(B, A).$ Further, by $(A\otimes_k V)\otimes_\m (M\otimes_k B)\cong M\otimes_kV\in \mathscr C,$ one gets
$M\otimes_k B\in {\rm mon}(B, \mathscr C).$

\vskip5pt

Conversely, if $M\otimes_k B\in {\rm mon}(B, \mathscr C),$ then by definition
$(A\otimes_kB)\otimes_\m(M\otimes_k B) \cong M\otimes_k B\in \mathscr C.$ Since $\mathscr C$ is closed under direct summands, it follows that $M\in \mathscr C.$

\vskip5pt

${\rm (ii)}$ \ Assume that $M\in \ ^\bot A$  and $M\notin \ {\rm gp}(A)$.
By (i), $M\otimes_kB\in {\rm mon}(B, \ ^\bot A) \subseteq  \ ^\bot \m$, where the inclusion follows from Theorem \ref{homodescriptionofsmon}.
Again by (i),  $M\otimes_kB\notin {\rm mon}(B, \ {\rm gp}(A)) =  {\rm gp}(\m)$, where the equality follows from Theorem \ref{gp}.

\vskip5pt

${\rm (iii)}$ \ Assume that $\mathscr C_1$ and $\mathscr C_2$  are additive full subcategories of $A$-mod closed under direct summands,
such that ${\rm mon}(B, \mathscr C_1) = {\rm mon}(B, \mathscr C_2)$.
We need to prove $\mathscr C_1 = \mathscr C_2$. Let $M\in \mathscr C_1$. By ${\rm (i)}$,  $M\otimes_k B\in {\rm mon}(B, \mathscr C_1).$
Thus $M\otimes_k B\in {\rm mon}(B, \mathscr C_2).$ Again by ${\rm (i)}$, $M\in \mathscr C_2.$ This completes the proof. $\s$

\vskip5pt

\begin{prop}\label{wgandmon2} \ If $\m$ is weakly Gorenstein, then so is $A$. Conversely, if $A$ is weakly Gorenstein, then a semi-Gorenstein-projective $\m$-module is Gorenstein-projective if and only if it is monic.

\vskip5pt

Thus, if $A$ is weakly Gorenstein, then $\m$ is weakly Gorenstein if and only if each semi-Gorenstein-projective $\m$-module is monic,
or equivalently, \ $^\bot \m = {\rm mon}(B, \ ^\perp A)$.
\end{prop}

\noindent{\bf Proof.} \ If $\m$ is weakly Gorenstein,  then $A$ is weakly Gorenstein, by  Lemma \ref{monmap}(ii).

\vskip5pt

Assume that $A$ is weakly Gorenstein and $M$ is a semi-Gorenstein-projective $\m$-module.
If $M$ is Gorenstein-projective, then $M$ is monic, by Theorem \ref{gp}. If $M$ is monic, then by Theorem \ref{homodescriptionofsmon} and Theorem \ref{gp} one has
\ $M \in {\rm mon}(B, A)\cap \ ^\bot \m = {\rm mon}(B, \ ^\bot A) = {\rm mon}(B, \ {\rm gp}(A)) = {\rm gp}(\m)$. $\s$

\vskip5pt

\begin{prop}\label{quiverwg1} \  Let $Q$ be a finite acyclic quiver. Then $A\otimes_kkQ$
is weakly Gorenstein if and only if $A$ is weakly Gorenstein.

\vskip5pt

In particular, $T_n(A) = \left(\begin{smallmatrix} A&A&\cdots &A&A \\ 0&A&\cdots &A&A \\  \vdots &\vdots & \vdots &\vdots & \vdots \\0& 0 & \cdots & 0 & A\end{smallmatrix}\right)$
is weakly Gorenstein if and only if $A$ is weakly Gorenstein. \end{prop}

\noindent{\bf Proof.} \ By Proposition \ref{wgandmon2}, it remains to prove the ``if" part. Assume that $A$ is weakly Gorenstein. We will prove that
$\m = A\otimes_kkQ$ is  weakly Gorenstein, by using induction on $|Q_0|$.
If $|Q_0| = 1$, then $\m = A$ is  weakly Gorenstein, by the assumption.

\vskip5pt

Assume that $|Q_0|\ge 2$. We write the conjunction of paths of $Q$ from left to right.
Since $Q$ is an acyclic quiver, $Q$ has a source vertex, say, $n$, and then $$kQ = \left(\begin{smallmatrix}
kQ'&{\rm rad}P(n)\\ 0&k \end{smallmatrix}\right)$$  where $Q'$ is the subquiver of $Q$ by deleting the source vertex $n$,
and $P(n) = kQe_n$. Then ${\rm rad}P(n)$ is a $kQ'$-$k$-bimodule. Thus $$\m = A\otimes_kkQ = \left(\begin{smallmatrix}
A\otimes_kkQ'&M\\ 0&A \end{smallmatrix}\right)$$
where $M = A\otimes_k{\rm rad}P(n)$ is an $(A\otimes_kkQ')$-$A$-bimodule. Since $Q$ is acyclic, so is $Q'$. Hence ${\rm rad}P(n)$ is a projective left $kQ'$-module. Thus
$M =  A\otimes_k {\rm rad}P(n)$ is a projective left $(A\otimes_kkQ')$-module, and also $M = A\otimes_k {\rm rad}P(n)$ is a projective right $A$-module.

\vskip5pt

Since $|Q'_0| = |Q_0|-1$, by induction $A\otimes kQ'$ is weakly Gorenstein. Applying Proposition \ref{matrixwg1} to
$\left(\begin{smallmatrix}
A\otimes kQ'&M\\ 0&A \end{smallmatrix}\right)= \m$, one sees that $\m$ is weakly Gorenstein. $\s$

\vskip5pt

\subsection{Semi-Groenstein-projective-free algebras}  Replacing the condition that $M_B$ is projective in Proposition \ref{matrixwg1} by ``$^\perp B = \add(B)$ and ${\rm proj.dim} M_B < \infty$", we then get the following
result on lsgp-free algebras.

\begin{prop}\label{matrixwg2} \  Assume that \ $^\perp B = \add(B)$ with ${\rm proj.dim} M_B < \infty$, and that $_AM$ is torsionless with \ ${\rm proj.dim} \ _AM < \infty$.
Then $\m=\left(\begin{smallmatrix}A&M\\0&B\end{smallmatrix}\right)$ is left weakly Gorenstein if and only if $A$ is left weakly Gorenstein.

\vskip5pt

Moreover, if in addition $_AM$ is projective, then \ $$^\perp \m = {\rm gp}(\m) = \{\binom{M\otimes_BP}{P}_{{\rm Id}_{M\otimes_BP}} \oplus \binom{G}{0} \ | \ P\in \add(B), \ G\in \ ^\perp A = {\rm gp}(A)\}$$ and \ $^\perp \m = \add(\m)$ if and only if \ $^\perp A = \add(A)$.\end{prop}

\noindent{\bf Proof.} \ Since \ $^\perp B = \add(B)$, it follows that \ $(^\perp B)^\perp   = B$-mod, and hence $\D(M_B)\in (^\perp B)^\perp$.
So, the conditions of  Proposition \ref{wgandmon1} are satisfied.
To prove the first assertion, by Proposition \ref{wgandmon1}, it suffices to prove that if $A$ is left weakly Gorenstein, then any semi-Gorenstein-projective $\m$-module $\binom{X}{Y}_\varphi$ is monic respect to bimodule \ $_AM_B$.
In fact, applying Theorem \ref{sgpovermatrixalg} to $\binom{X}{Y}_\varphi\in \ ^\perp \m$ one gets that \ $Y\in \ ^\perp B$, that $\varphi: M\otimes_BY\longrightarrow X$
induces isomorphisms $\Ext^i_A(X, A)\cong\Ext^i_A(M\otimes_BY, A)$ for all $i\geq1$, and that $\varphi^*: X^*\longrightarrow (M\otimes_BY)^*$ is a right $A$-epimorphism.
Since $Y\in \ ^\perp B = \add(B)$, \ $_BY$ is projective. Thus \ $_A(M\otimes_BY)\in \add(_AM)$. Since by assumption $_AM$ is torsionless,
it follows that \ $_A(M\otimes_BY)$ is torsionless. Thus, the canonical map
$\phi_{M\otimes_BY}: M\otimes_BY \longrightarrow (M\otimes_BY)^{**}$ is a monomorphism.
By the same argument as in the proof of Proposition \ref{matrixwg1} one concludes that $\binom{X}{Y}_\varphi$ is monic respect to bimodule \ $_AM_B$.

\vskip10pt

Now, assume in addition that $_AM$ is projective. Continuing the argument above, one knows that \ $_A(M\otimes_BY)$ is projective, thus,
$\phi_{M\otimes_BY}: M\otimes_BY \longrightarrow (M\otimes_BY)^{**}$ is an isomorphism.
By $\Ext^i_A(X, A)\cong\Ext^i_A(M\otimes_BY, A) = 0$ for all $i\geq1$, one has $X\in \ ^\perp A$. Since by
assumption $A$ is weakly Gorenstein, $X\in {\rm gp}(A)$, and hence $\phi_X: X\longrightarrow X^{**}$ is an isomorphism.
Since $\varphi^*: \ X^* \longrightarrow (M\otimes_B Y)^*$ is an epimorphism and $(M\otimes_B Y)^*$ is a right projective $A$-module, it follows that $\varphi^*$ is a splitting epimorphism, and hence
$\varphi^{**}: \ (M\otimes_B Y)^{**} \longrightarrow X^{**}$ is a splitting monomorphism. From the commutative diagram
$$\xymatrix{
M\otimes_B Y\ar[d]_{\phi_{M\otimes_B Y}}^-{\cong} \ar[r]^-{\varphi} & X \ar[d]^-{\phi_{X}}_-{\cong} \\
(M\otimes_B Y)^{**} \ar@{^{(}->}[r]^-{\varphi^{**}} & X^{**}}$$
one sees that $\varphi: M\otimes_B Y \longrightarrow X$ is also a splitting monomorphism. Thus $X \cong (M\otimes_BY) \oplus X'$ for some $X'\in {\rm gp}(A)$ and
$\binom{X}{Y}_\varphi = \binom{M\otimes_BY}{Y}_{{\rm Id}_{M\otimes_BY}} \oplus \binom{X'}{0}$.

\vskip5pt

Since \ $_AM$ is projective and ${\rm proj.dim}M_B <\infty$, the $A$-$B$-bimodule $M$ is compatible (cf. Lemma \ref{compatible}).
By Theorem \ref{gpovermatrixalg},  $\binom{X'}{0}_\varphi \in {\rm gp}(\m)$, and hence
$\binom{X}{Y}_\varphi = \binom{M\otimes_BY}{Y}_{{\rm Id}_{M\otimes_BY}} \oplus \binom{X'}{0}\in {\rm gp}(\m)$.
This proves $$^\perp\m ={\rm gp}(\m) = \{\binom{M\otimes_BP}{P}_{{\rm Id}_{M\otimes_BP}} \oplus \binom{G}{0} \ | \ P\in \add(B), \ G\in \ ^\perp A = {\rm gp}(A)\}$$
and from which one sees that $^\perp \m = \add(\m)$ if and only if \ $^\perp A = \add(A)$. $\s$

\begin{rem} \ The ``only if" part in {\rm Proposition \ref{matrixwg2}} does not need the conditions that \ $^\perp B = \add(B)$ and $_AM$ is torsionless.
\end{rem}

\begin{prop}\label{quiverwg2} \  Let $Q$ be a finite acyclic quiver, $I$ an admissible ideal of $kQ$, and $\m = A\otimes_kkQ/I$.
Then \ $^\perp \m = \add(\m)$ if and only if \ $^\perp A = \add(A)$.
\end{prop}

\noindent{\bf Proof.} \  Assume that \ $^\perp A = \add(A)$. We will prove \ $^\perp\m = \add(\m)$,
again by using induction on $|Q_0|$.
If $|Q_0| = 1$, then $\m = A$, thus the assertion holds, by the assumption
$^\perp A = \add(A)$.

\vskip5pt

Assume that $|Q_0|\ge 2$. Similar as in the proof of Proposition \ref{quiverwg1}, we write $\m$ as a triangular matrix algebra.
However, in order to apply Proposition \ref{matrixwg2}, this time we need to use
the subquiver $Q'$ of $Q$ by deleting a sink vertex, say, $1$,  and the corresponding algebra $kQ'/I'$.
Then $$kQ/I = \left(\begin{smallmatrix}
k&{\rm rad}(e_1kQ/I) \\ 0& kQ'/I' \end{smallmatrix}\right)$$ where ${\rm rad}(e_1kQ/I)$ is a $k$-$(kQ'/I')$-bimodule. Thus $$\m = A\otimes_k(kQ/I) = \left(\begin{smallmatrix}
A&M\\ 0&\m' \end{smallmatrix}\right)$$
where $\m' = A\otimes_k(kQ'/I'), \ M = A\otimes_k{\rm rad}(e_1kQ/I)$ is a $A$-$\m'$-bimodule.
Since $|Q_0| = |Q_0| - 1$, by induction one gets \ $^\perp \m' = \add(\m').$

\vskip5pt

Since ${\rm proj.dim.}{\rm rad}(e_1kQ/I)_{kQ'/I'} <\infty$, it follows that
$${\rm proj.dim.} M_{\m'} < \infty.$$
Also, $_AM$ is projective. Since we already known \ $^\perp \m' = \add(\m')$ by induction, thus we can apply Proposition \ref{matrixwg2} to  $\m = A\otimes_k(kQ/I) = \left(\begin{smallmatrix}
A&M\\ 0&\m' \end{smallmatrix}\right)$ to get
$$^\perp \m = {\rm gp}(\m) = \{\binom{M\otimes_{\m'}P}{P}_{{\rm Id}_{M\otimes_{\m'}P}} \oplus \binom{G}{0} \ | \ P\in \add(\m'), \ G\in \ ^\perp A\}.$$
Since $G\in \ ^\perp A = \add(A)$, it follows that $\binom{G}{0}\in \add(\m)$, and hence \ $^\perp \m  = \add(\m)$.

\vskip5pt

Conversely, assume that \ $^\perp\m = \add(\m)$. Let $X$ be an indecomposable $A$-module with $X\in \ ^\perp A$.
For any indecomposable projective $(kQ/I)$-module $P$, by the Cartan-Eilenberg isomorphism
([CE, Thm. 3.1, p.209, p.205]) one has
$${\rm Ext}^{i}_{\m}(X\otimes_k
P, A\otimes_k kQ/I) = \bigoplus\limits_{p+q = i}({\rm Ext}^p_A(X,
A)\otimes_k {\rm Ext}^q_{kQ/I}(P, kQ/I)) = 0, \ \forall \ i\ge 1.$$
So $X\otimes_kP\in \ ^\perp \m = \add(\m)$, and hence $X\in \add(A).$ $\s$

\subsection{Example and Problem} \ Let $A$
be the algebra given by quiver $\xymatrix{\ar@(ur, ul)^\beta 2
\ar@{->}@<2pt> [r]^-{\alpha} &*-<1pt> {\begin{matrix} \ 1
  \end{matrix}}}$ and relations $\beta^2, \ \alpha\beta$. The Auslander-Reiten
quiver of $A$ is
$$
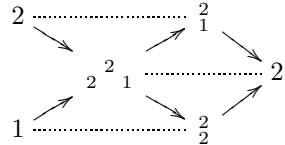
\xymatrix@C=1.2pc @ R = 0.3pc{2\ar[dr]\ar@{.}[rr]&& {\begin{smallmatrix}2\\1\end{smallmatrix}}\ar[dr]
\\ & {\begin{smallmatrix}&2\\2&&1\end{smallmatrix}}\ar[dr]\ar[ur]\ar@{.}[rr] & & 2
\\ 1 \ar[ur] \ar@{.}[rr]&& {\begin{smallmatrix}2\\2\end{smallmatrix}}\ar[ur]}$$ with indecomposable projective modules $P(1) = 1$ and $P(2) =
\begin{smallmatrix}&2\\2&&1\end{smallmatrix}$, and indecomposable injective modules $I(1) = \begin{smallmatrix}2\\1\end{smallmatrix}$ and $I(2) =
\begin{smallmatrix}2\\2\end{smallmatrix}$.
Since $$\Ext^1_A(2, \
\begin{smallmatrix}&2\\2&&1\end{smallmatrix})\ne 0, \ \ \ \Ext^1_A(\begin{smallmatrix}2\\2\end{smallmatrix}, \
1) \ne 0, \ \ \ \Ext^2_A(\begin{smallmatrix}2\\1\end{smallmatrix}, \
1) = \Ext^1_A(\begin{smallmatrix}2\\1\end{smallmatrix}, \ 2)\ne 0$$
one sees that $A$ is lsgp-free, i.e., $^\perp A = \add(A)$.  Note that $A$ is not Gorenstein. By Proposition \ref{quiverwg2},
$^\perp (A\otimes_kkQ/I) = {\rm gp}(A\otimes_kkQ/I) = \add(A\otimes_kkQ/I)$, for any finite acyclic quiver $Q$ and any admissible ideal $I$.

\vskip10pt

{\bf Problem 1.} \ Are there a left weakly Gorenstein algebra $A$, a finite acyclic quiver $Q$, and an admissible ideal $I$ of $kQ$,
such that $A\otimes_kkQ/I$ is not left weakly Gorenstein, or equivalently, such that there is a semi-Gorenstein-projective $(A\otimes_kkQ/I)$-module which is not monic?

\vskip5pt

This  problem is different from {\bf Question 4.} Such an algebra $A$ (if there exists) is not Gorenstein (otherwise, $A\otimes_k kQ/I$ is Gorenstein, by [AR, Proposition 2.2]); such an $I\ne 0$, by Proposition \ref{quiverwg1};
and also \ $^\perp A \ne \add(A)$, by Proposition \ref{quiverwg2}.

\vskip5pt

\section{\bf Canonical maps of modules  over $T_2(A)$}

Let $A$ be an Artin algebra, $\m = T_2(A) = \left(\begin{smallmatrix}A&A\\0&A\end{smallmatrix}\right),$ and $M$ a $\m$-module.
We will give a sufficient and necessary condition, such that the canonical
$\m$-map $\phi_M: M \longrightarrow M^{**}$ is a monomorphism (an epimorphism, and reflexive, respectively);
and we will give a sufficient and necessary condition such that $M$ is double semi-Gorenstein-projective with $\phi_M$ a monomorphism (an epimorphism, respectively).

\vskip5pt

Recall that a left $\m$-module $M$ is identified with the triple $\left(\begin{smallmatrix}X\\ Y\end{smallmatrix}\right)_\varphi$, where $\varphi: Y \longrightarrow X$ is a left $A$-map;
and a right $\m$-module is identified with a triple $(U, V)_\psi$, where $\psi: U\longrightarrow V$ is a right $A$-map.
Using the identifications, we will determine the right $\m$-module $M^* = \Hom_\m(M, \ _\m\m)$, the left $\m$-module $M^{**} = \Hom_\m(M^*, \ \m_\m)$, and  $\phi_M: M \longrightarrow M^{**}$.

\subsection{The $\m$-dual of a left $\m$-module}  For a left $\m$-module $\left(\begin{smallmatrix}X\\ Y\end{smallmatrix}\right)_\varphi$,
we will determine the right $\m$-module $\left(\begin{smallmatrix}X\\ Y\end{smallmatrix}\right)^*_\varphi
= \Hom_\m(\left(\begin{smallmatrix}X\\ Y\end{smallmatrix}\right)_\varphi, \ _\m\m)$.
As a left $\m$-module, $_\m\m = \left(\begin{smallmatrix}A\\ 0\end{smallmatrix}\right)\oplus \left(\begin{smallmatrix}A\\ A\end{smallmatrix}\right)_{{\rm Id}_A} =
\left(\begin{smallmatrix}A\oplus A\\ A\end{smallmatrix}\right)_{\binom{0}{{\rm Id}_A}}$. Thus, any $\m$-map $$f\in \left(\begin{smallmatrix}X\\ Y\end{smallmatrix}\right)^*_\varphi = \Hom_\m(\left(\begin{smallmatrix}X\\ Y\end{smallmatrix}\right)_\varphi, \left(\begin{smallmatrix}A\oplus A\\ A\end{smallmatrix}\right)_{\binom{0}{{\rm Id}_A}})$$
 is of the form
$\left(\begin{smallmatrix}\left(\begin{smallmatrix}\alpha_1\\ \alpha_2\end{smallmatrix}\right)\\ \beta\end{smallmatrix}\right)$, where $\alpha_1, \alpha_2\in X^* = \Hom_A(X, A), \ \beta\in Y^*$,
such that the square $$\xymatrix{Y\ar[d]_{\beta} \ar[r]^-{\varphi} & X \ar[d]^-{\binom{\alpha_1}{\alpha_2}} \\
A\ar[r]^-{\binom{0}{1}} & A\oplus A}$$
commutes. So $\alpha_1\varphi = 0$, $\beta = \alpha_2\varphi.$ Thus, there is a unique $g\in ({\rm Coker} \varphi)^* = \Hom_A({\rm Coker} \varphi, A)$ such that
$\alpha_1 = g\pi,$ where $\pi: X\longrightarrow {\rm Coker}\varphi$ is the canonical $A$-epimorphism.

\vskip5pt

\begin{lem} \label{T2dualofleftmod} \ Let $\left(\begin{smallmatrix}X\\ Y\end{smallmatrix}\right)_\varphi$ be a left $\m$-module with $\varphi: Y \longrightarrow X$ a left $A$-map. Then
\vskip5pt

$(\rm i)$ \ Any $f\in \left(\begin{smallmatrix}X\\ Y\end{smallmatrix}\right)^*_\varphi$ is of the form
$\left(\begin{smallmatrix}\left(\begin{smallmatrix}g\pi\\ \alpha_2\end{smallmatrix}\right)\\ \alpha_2\varphi\end{smallmatrix}\right)$, where $g\in ({\rm Coker} \varphi)^*$,
$\pi: X\longrightarrow {\rm Coker}\varphi$ is the canonical $A$-epimorphism,  and $\alpha_2\in X^*$.

\vskip5pt

$(\rm ii)$ \ There is a unique right $\m$-module isomorphism \ \
$h: \left(\begin{smallmatrix}X\\ Y\end{smallmatrix}\right)_\varphi^* \cong (({\rm Coker} \varphi)^*, \ X^*)_{\pi^*}$, given by
$$f= \left(\begin{smallmatrix}\left(\begin{smallmatrix}g\pi\\ \alpha_2\end{smallmatrix}\right)\\ \alpha_2\varphi\end{smallmatrix}\right)\mapsto (g, \alpha_2)$$
where $\pi^*: ({\rm Coker} \varphi)^*\longrightarrow X^*$ is the right $A$-monomorphism induced by $\pi$.
\end{lem}

\noindent{\bf Proof.} $(\rm ii)$ \ We claim that $h$ is a right $\m$-map, i.e.,
$$h(f\left(\begin{smallmatrix}a_1 & a_2\\ 0 & a_3\end{smallmatrix}\right)) = h(f)\left(\begin{smallmatrix}a_1 & a_2\\ 0 & a_3\end{smallmatrix}\right), \ \ \forall \ f = \left(\begin{smallmatrix}\left(\begin{smallmatrix}g\pi\\ \alpha_2\end{smallmatrix}\right)\\ \alpha_2\varphi\end{smallmatrix}\right)\in \left(\begin{smallmatrix}X\\ Y\end{smallmatrix}\right)_\varphi^*, \ \ \forall \ \left(\begin{smallmatrix}a_1 & a_2\\ 0 & a_3\end{smallmatrix}\right)\in \m.$$ In fact, for any $\binom{x}{y}\in \binom{X}{Y}_\varphi$, since $\pi \varphi = 0$, one has
\begin{align*} (f\left(\begin{smallmatrix}a_1 & a_2\\ 0 & a_3\end{smallmatrix}\right))\binom{x}{y}
& = (f\binom{x}{y})\left(\begin{smallmatrix}a_1 & a_2\\ 0 & a_3\end{smallmatrix}\right)
= \left(\begin{smallmatrix}(g\pi)(x) & \alpha_2(x)\\ 0 & \alpha_2(\varphi(y))\end{smallmatrix}\right)\left(\begin{smallmatrix}a_1 & a_2\\ 0 & a_3\end{smallmatrix}\right)
\\ & = \left(\begin{smallmatrix}(g\pi)(x)a_1 & (g\pi)(x)a_2 + \alpha_2(x)a_3\\ 0 & \alpha_2(\varphi(y))a_3\end{smallmatrix}\right)
= \left(\begin{smallmatrix}((ga_1)\pi)(x) & ((ga_2)\pi + \alpha_2a_3)(x)\\ 0 & ((ga_2)\pi  + \alpha_2a_3)(\varphi(y))\end{smallmatrix}\right)
\\ & = \left(\begin{smallmatrix}\left(\begin{smallmatrix}(ga_1)\pi \\ (ga_2)\pi  + \alpha_2a_3\end{smallmatrix}\right)\\ ((ga_2)\pi  + \alpha_2a_3)\varphi\end{smallmatrix}\right)\binom{x}{y}
\end{align*}
Thus $f\left(\begin{smallmatrix}a_1 & a_2\\ 0 & a_3\end{smallmatrix}\right) = \left(\begin{smallmatrix}\left(\begin{smallmatrix}(ga_1)\pi \\ (ga_2)\pi  + \alpha_2a_3\end{smallmatrix}\right)\\ ((ga_2)\pi  + \alpha_2a_3)\varphi\end{smallmatrix}\right),$ and hence $$h(f\left(\begin{smallmatrix}a_1 & a_2\\ 0 & a_3\end{smallmatrix}\right)) = (ga_1, (ga_2)\pi + \alpha_2a_3).$$
One the other hand, by the right $\m$-module structure of  $(({\rm Coker} \varphi)^*, X^*)_{\pi^*}$, one has
\begin{align*} h(f)\left(\begin{smallmatrix}a_1 & a_2\\ 0 & a_3\end{smallmatrix}\right) & = (g, \alpha_2)\left(\begin{smallmatrix}a_1 & a_2\\ 0 & a_3\end{smallmatrix}\right)
=(ga_1, \pi^*(ga_2) + \alpha_2a_3) \\ & = (ga_1, (ga_2)\pi + \alpha_2a_3).\end{align*}
This proves the claim.

\vskip5pt

Since the map
$$(({\rm Coker} \varphi)^*, X^*)_{\pi^*}\longrightarrow \left(\begin{smallmatrix}X\\ Y\end{smallmatrix}\right)_\varphi^*, \ \ \
(g, \alpha_2)\mapsto f= \left(\begin{smallmatrix}\left(\begin{smallmatrix}g\pi\\ \alpha_2\end{smallmatrix}\right)\\ \alpha_2\varphi\end{smallmatrix}\right)$$
is the inverse of $h$, $h$ is a right $\m$-isomorphism. $\s$

\subsection{The $\m$-dual of a right $\m$-module}  \ Similarly, one can determine the
$\m$-dual of a right $\m$-module
$(U, V)_\psi$, where $\psi: U \longrightarrow V$ is a right $A$-map. As a right $\m$-module, $\m_\m = (A, A)_{{\rm Id}_A}\oplus (0, A) = (A, A\oplus A)_{\binom{{\rm Id}_A}{0}}$. So any $\m$-map
$$f\in (U, V)_\psi^* = \Hom_\m((U, V)_\psi, (A, A\oplus A)_{\binom{{\rm Id}_A}{0}})$$
is of the form $(\alpha, \left(\begin{smallmatrix}\beta_1\\ \beta_2\end{smallmatrix}\right))$, where $\alpha\in U^* = \Hom_A(U, A_A), \ \beta_1, \ \beta_2\in V^*$,
such that $$\xymatrix{U\ar[d]_{\alpha} \ar[r]^-{\psi} & V \ar[d]^-{\binom{\beta_1}{\beta_2}} \\
A\ar[r]^-{\binom{1}{0}} & A\oplus A}$$
commutes. Thus $\alpha = \beta_1\psi, \ \beta_2\psi = 0.$ Hence, there is a unique $g\in ({\rm Coker} \psi)^* = \Hom_A({\rm Coker} \psi, \ A_A)$ such that
$\beta_2 = g\pi,$ where $\pi: V\longrightarrow {\rm Coker}\psi$ is the canonical $A$-map.
By the similar argument  one has

\vskip5pt

\begin{lem} \label{T2dualofrightmod} \ Let $(U, V)_\psi$ be a right $\m$-module with $\psi: U \longrightarrow V$ a right $A$-map. Then

\vskip5pt

$(\rm i)$ \ Any $f\in (U, V)_\psi^*$ is of the form
$(\beta_1\psi, \left(\begin{smallmatrix}\beta_1\\ g\pi\end{smallmatrix}\right))$, where $\beta_1\in V^*$, $g\in ({\rm Coker} \psi)^*$,
and $\pi: V\longrightarrow {\rm Coker}\psi$ is the canonical $A$-epimorphism.

\vskip5pt

$(\rm ii)$ \ There is a unique left $\m$-module isomorphism \ \ $h': (U, V)_\psi^* \cong \left(\begin{smallmatrix}V^* \\ ({\rm Coker} \psi)^*\end{smallmatrix}\right)_{\pi^*}$, given by
$$f= (\beta_1\psi, \left(\begin{smallmatrix}\beta_1\\ g\pi\end{smallmatrix}\right))\mapsto \binom{\beta_1}{g}$$
where $\pi^*: ({\rm Coker} \psi)^*\longrightarrow V^*$ is the left $A$-monomorphism induced by $\pi$.
\end{lem}

\subsection{The left $\m$-module $\binom{X}{Y}^{**}_{\varphi}$} \ For any left $\m$-module $\left(\begin{smallmatrix}X\\ Y\end{smallmatrix}\right)_\varphi$
with left $A$-map $\varphi: Y \longrightarrow X$, by Lemma \ref{T2dualofleftmod}, one has the right module isomorphism
$$h: \binom{X}{Y}^{*}_{\varphi} \cong (({\rm Coker} \varphi)^*, \ X^*)_{\pi^*}, \ \ \
f= \left(\begin{smallmatrix}\left(\begin{smallmatrix}g\pi\\ \alpha_2\end{smallmatrix}\right)\\ \alpha_2\varphi\end{smallmatrix}\right)\mapsto (g, \alpha_2)\eqno(5.1)$$
where $\pi^*: ({\rm Coker} \varphi)^*\longrightarrow X^*$ is the right $A$-monomorphism induced by $\pi: X\longrightarrow {\rm Coker}\varphi$.
Applying Lemma \ref{T2dualofrightmod} to $(({\rm Coker} \varphi)^*, \ X^*)_{\pi^*}$, we then get

\vskip5pt

\begin{lem} \label{doubledual} \ $(\rm i)$ \ Any $f\in (({\rm Coker} \varphi)^*, \ X^*)^*_{\pi^*}$ is of the form
$(\beta_1\pi^*, \left(\begin{smallmatrix}\beta_1\\ gp\end{smallmatrix}\right))$, where $\beta_1\in X^{**}$, \ $g\in ({\rm Coker} \pi^*)^*$,
and $p: X^*\longrightarrow {\rm Coker} \pi^*$ is the canonical $A$-epimorphism.

\vskip5pt

$(\rm ii)$ \ There is a unique left $\m$-module isomorphism \ \ $\widetilde{h}: \ \left(\begin{smallmatrix}X^{**} \\ ({\rm Coker} \pi^*)^*\end{smallmatrix}\right)_{p^*}\cong \left(\begin{smallmatrix}X\\ Y\end{smallmatrix}\right)_\varphi^{**}$, \ given by
$$\binom{\beta_1}{g}\mapsto h^*((\beta_1\pi^*, \left(\begin{smallmatrix}\beta_1\\ gp\end{smallmatrix}\right)))$$
where $p^{*}: ({\rm Coker} \pi^*)^*\longrightarrow X^{**}$ is the $A$-monomorphism induced by $p$,  \ \ $h: \binom{X}{Y}^{*}_{\varphi} \cong (({\rm Coker} \varphi)^*, \ X^*)_{\pi^*}$ is given in $(5.1)$, and \ \ $h^*: (({\rm Coker} \varphi)^*, \ X^*)_{\pi^*}^*\longrightarrow \binom{X}{Y}^{**}_{\varphi} $ is induced by $h$.
\end{lem}

\subsection{The canonical $\m$-map $\phi_{\binom{X}{Y}_\varphi}: \binom{X}{Y}_\varphi\longrightarrow \binom{X}{Y}_\varphi^{**}$} \ \ For a left $\m$-module $\left(\begin{smallmatrix}X\\ Y\end{smallmatrix}\right)_\varphi$, one has an exact sequence $Y \stackrel \varphi \longrightarrow X \stackrel \pi \longrightarrow {\rm Coker}\varphi \longrightarrow 0$ of left $A$-modules.
Applying $\Hom_A(-, \ _AA)$, one gets an exact sequence of right $A$-modules
$$0 \longrightarrow ({\rm Coker}\varphi)^* \stackrel {\pi^*} \longrightarrow X^* \stackrel {\varphi^*} \longrightarrow Y^*$$
and the exact sequence
$$0 \longrightarrow ({\rm Coker}\varphi)^* \stackrel {\pi^*} \longrightarrow X^* \stackrel {p} \longrightarrow {\rm Coker}\pi^* \longrightarrow 0.$$
Thus, there is a unique $A$-map $\beta: {\rm Coker}\pi^* \longrightarrow Y^*$ such that the diagram
$$\xymatrix{0\ar[r] & ({\rm Coker}\varphi)^*\ar[r]^-{\pi^*}\ar@{=}[d] &  X^{*} \ar[r]^-{p}\ar@{=}[d] & {\rm Coker}\pi^*\ar@{-->}[d]^{\beta} \ar[r] & 0  \\
0 \ar[r] & ({\rm Coker}\varphi)^*\ar[r]^-{\pi^*} & X^* \ar[r]^-{\varphi^*} & Y^*.}\eqno(5.2)$$
commutes, i.e., $\varphi^* = \beta p.$ Thus, $\varphi^*$ is an epimorphism if and only if so is $\beta$, and if and only if $\beta$ is an isomorphism.
So one has the $A$-map $\beta^*: Y^{**}\longrightarrow ({\rm Coker}\pi^*)^*.$ Consider the composition
$$\beta^*\phi_Y: Y\longrightarrow ({\rm Coker}\pi^*)^*$$
where $\phi_Y: Y\longrightarrow Y^{**}$ is the canonical map. By the definition of $\phi_Y$ and $\beta^*$, one knows that
$\beta^*\phi_Y: Y \longrightarrow ({\rm Coker}\pi^*)^*$ is given by
$$y\mapsto ``g\mapsto (\beta(g))(y)", \  \ \forall \ g\in {\rm Coker}\pi^*$$
i.e., $((\beta^*\phi_Y)(y))(g) = (\beta(g))(y), \ \forall \ y\in Y.$

\vskip5pt

\begin{prop} \label{thecanonicalmap} \ For any left $\m$-module $\left(\begin{smallmatrix}X\\ Y\end{smallmatrix}\right)_\varphi$
with left $A$-map $\varphi: Y \longrightarrow X$, with the notations above one has

\vskip5pt

$(\rm i)$ \ \ $\binom{\phi_X}{\beta^*\phi_Y}: \binom{X}{Y}_\varphi \longrightarrow \binom{X^{**}}{({\rm Coker}\pi^*)^*}_{p^*}$ is left $\m$-map,
where $\phi_X: X \longrightarrow X^{**}$ and $\phi_Y: Y \longrightarrow Y^{**}$ are the canonical $A$-maps,
$\beta: {\rm Coker}\pi^* \longrightarrow Y^*$ is the canonical $A$-map such that $\varphi^* = \beta p$, and $\beta^*: Y^{**}\longrightarrow ({\rm Coker}\pi^*)^*$ is induced by $\beta$.

\vskip5pt

$(\rm ii)$ \ The canonical $\m$-map $\phi_{\binom{X}{Y}_\varphi}: \binom{X}{Y}_\varphi\longrightarrow \binom{X}{Y}_\varphi^{**}$ is given by
$$\phi_{\binom{X}{Y}_\varphi} = \widetilde{h}\circ \binom{\phi_X}{\beta^*\phi_Y}.$$
where $\widetilde{h}: \ \left(\begin{smallmatrix}X^{**} \\ ({\rm Coker} \pi^*)^*\end{smallmatrix}\right)_{p^*} \longrightarrow \binom{X}{Y}_\varphi^{**}$ is the isomorphism given in {\rm Lemma  \ref {doubledual}}.
\end{prop}
\noindent{\bf Proof.} \ $(\rm i)$ \ One needs to prove the diagram
$$\xymatrix{Y\ar[r]^-{\varphi}\ar[d]_-{\beta^*\phi_Y} &  X\ar[d]^-{\phi_X}\\
({\rm Coker}\pi^*)^*\ar[r]^-{p^*} & X^{**}}$$
commutes, i.e., $p^*\beta^*\phi_Y = \phi_X \varphi$. In fact, since $\varphi^* = \beta p$, one has $\varphi^{**} = p^*\beta^*.$
By the functorial property of the canonical map $\phi_X: X \longrightarrow X^{**}$ one has the commutative diagram
$$\xymatrix{Y\ar[r]^-{\varphi}\ar[d]_-{\phi_Y} &  X\ar[d]^-{\phi_X}\\
Y^{**}\ar[r]^-{\varphi^{**}} & X^{**}.}$$
It follows that \ $p^*\beta^*\phi_Y = \varphi^{**}\phi_Y = \phi_X \varphi.$

\vskip5pt

$(\rm ii)$ \ We need to prove $\phi_{\binom{X}{Y}_\varphi}\binom{x}{y} = \widetilde{h}(\binom{\phi_X}{\beta^*\phi_Y}\binom{x}{y}), \ \ \forall \ \binom{x}{y}\in \binom{X}{Y}_\varphi.$
\ For this, let \ $f\in \binom{X}{Y}_\varphi^*$. By Lemma \ref{T2dualofleftmod}(i), \ $f = \left(\begin{smallmatrix}\left(\begin{smallmatrix}g\pi\\ \alpha_2\end{smallmatrix}\right)\\ \alpha_2\varphi\end{smallmatrix}\right)$, where \ $g\in ({\rm Coker} \varphi)^*$, \
$\pi: X\longrightarrow {\rm Coker}\varphi$ \ is the canonical $A$-epimorphism,  and \ $\alpha_2\in X^*$. By the definition of $\phi_{\binom{X}{Y}_\varphi}$ one has
$$\phi_{\binom{X}{Y}_\varphi}\binom{x}{y}(f) = f\binom{x}{y} = \left(\begin{smallmatrix}\left(\begin{smallmatrix}g\pi\\ \alpha_2\end{smallmatrix}\right)\\ \alpha_2\varphi\end{smallmatrix}\right)\binom{x}{y}
= \left(\begin{smallmatrix}(g\pi)(x) & \alpha_2(x)\\ 0& \alpha_2(\varphi(y))\end{smallmatrix}\right)\in \m.$$

\vskip5pt
\noindent On the other hand,  by the definitions of $\beta^*\phi_Y$ and $\widetilde{h}$ one has
\begin{align*}&\widetilde{h}(\binom{\phi_X}{\beta^*\phi_Y}\binom{x}{y})(f) =\widetilde{h}(\binom{\phi_X(x)}{\beta^*\phi_Y(y)})(f)
\\ & = h^*((\phi_X(x)\pi^*, \left(\begin{smallmatrix}\phi_X(x)\\ (\beta^*\phi_Y(y)) p\end{smallmatrix}\right)))(f)
= (\phi_X(x)\pi^*, \left(\begin{smallmatrix}\phi_X(x)\\ (\beta^*\phi_Y(y)) p\end{smallmatrix}\right))(h(f))
\\ &= (\phi_X(x)\pi^*, \left(\begin{smallmatrix}\phi_X(x)\\ (\beta^*\phi_Y(y)) p\end{smallmatrix}\right)) (g, \alpha_2)
= \left(\begin{smallmatrix}\phi_X(x)(\pi^*(g)) & \phi_X(x)(\alpha_2)\\ 0& \beta^*\phi_Y(y)(p(\alpha_2))\end{smallmatrix}\right)
\\ & = \left(\begin{smallmatrix}\pi^*(g)(x) & \alpha_2(x)\\ 0& \phi_Y(y)(\beta(p(\alpha_2)))\end{smallmatrix}\right)
= \left(\begin{smallmatrix}\pi^*(g)(x) & \alpha_2(x)\\ 0& \beta(p(\alpha_2))(y)\end{smallmatrix}\right)
\\ & = \left(\begin{smallmatrix}g(\pi(x)) & \alpha_2(x)\\ 0& \varphi^*(\alpha_2)(y)\end{smallmatrix}\right)
= \left(\begin{smallmatrix}(g\pi)(x) & \alpha_2(x)\\ 0& \alpha_2(\varphi(y))\end{smallmatrix}\right)\in \m.\end{align*}
This
completes the proof. $\s$

\subsection{Torsionless $\m$-modules and reflexive $\m$-modules}

\vskip5pt

\begin{cor} \label{reflexive} \ Let  $\left(\begin{smallmatrix}X\\ Y\end{smallmatrix}\right)_\varphi$
be a left $\m$-module, where $\varphi: Y \longrightarrow X$ is a left $A$-map. Then

\vskip5pt

$(\rm i)$ \ \ $\left(\begin{smallmatrix}X\\ Y\end{smallmatrix}\right)_\varphi$
is a torsionless $\m$-module if and only if it is monic $($i.e., $\varphi$ is a monomorphism$)$,
$X$ and $Y$ are torsionless $A$-modules.

\vskip5pt

$(\rm ii)$ \ \ $\phi_{\left(\begin{smallmatrix}X\\ Y\end{smallmatrix}\right)_\varphi}$
is a $\m$-epimorphism if and only if $\phi_X$ and $\beta^*\phi_Y: Y \longrightarrow ({\rm Coker}\pi^*)^*$ are $A$-epimorphisms.

\vskip5pt

$(\rm iii)$ \ \ $\left(\begin{smallmatrix}X\\ Y\end{smallmatrix}\right)_\varphi$
is a reflexive $\m$-module if and only if $\varphi$ is a monomorphism,
$X$ is reflexive $A$-module, and $\beta^*\phi_Y$ is an isomorphism.
\end{cor}
\noindent{\bf Proof.} \ By Proposition \ref{thecanonicalmap}, $\phi_{\binom{X}{Y}_\varphi}: \binom{X}{Y}_\varphi\longrightarrow \binom{X}{Y}_\varphi^{**}$ is given by
$$\phi_{\binom{X}{Y}_\varphi} = \widetilde{h}\circ \binom{\phi_X}{\beta^*\phi_Y}$$
where $\widetilde{h}: \ \left(\begin{smallmatrix}X^{**} \\ ({\rm Coker} \pi^*)^*\end{smallmatrix}\right)_{p^*} \longrightarrow \binom{X}{Y}_\varphi^{**}$ is an isomorphism.
Since $\binom{\phi_X}{\beta^*\phi_Y}: \binom{X}{Y}_\varphi \longrightarrow \binom{X^{**}}{({\rm Coker}\pi^*)^*}_{p^*}$ is a left $\m$-map, the diagram
$$\xymatrix{Y\ar[r]^-{\varphi}\ar[d]_-{\beta^*\phi_Y} &  X\ar[d]^-{\phi_X}\\
({\rm Coker}\pi^*)^*\ar[r]^-{p^*} & X^{**}}$$
commutes, i.e., \ $p^*\beta^*\phi_Y = \phi_X \varphi$.

\vskip5pt

$(\rm i)$ \ Assume that $\left(\begin{smallmatrix}X\\ Y\end{smallmatrix}\right)_\varphi$
is a torsionless $\m$-module, i.e.,  $\phi_{\binom{X}{Y}_\varphi}: \binom{X}{Y}_\varphi\longrightarrow \binom{X}{Y}_\varphi^{**}$  is a $\m$-monomorphism.
Thus $\phi_X$ and $\beta^*\phi_Y$ are monomorphisms,
in particular  $\phi_Y$ is a monomorphism, so $X$ and $Y$ are torsionless.

Since $p^*\beta^*\phi_Y = \phi_X \varphi$ and $p^*: ({\rm Coker}\pi^*)^* \longrightarrow X^{**}$ is a monomorphism,
it follows that $\phi_X \varphi$ is a monomorphism, and hence $\varphi$ is a monomorphism, i.e.,
$\left(\begin{smallmatrix}X\\ Y\end{smallmatrix}\right)_\varphi$ is monic.

\vskip5pt

Conversely, assume that $\left(\begin{smallmatrix}X\\ Y\end{smallmatrix}\right)_\varphi$ is monic, $X$ and $Y$ are torsionless $A$-modules, i.e.,
$\varphi$ is a monomorphism, $\phi_X$ and $\phi_Y$ are monomorphisms.
By $p^*\beta^*\phi_Y = \phi_X \varphi$, \ $p^*\beta^*\phi_Y$ is a monomorphism, and hence $\beta^*\phi_Y$ is a monomorphism. Thus $\phi_{\binom{X}{Y}_\varphi} = \widetilde{h}\circ \binom{\phi_X}{\beta^*\phi_Y}$ is a monomorphism, i.e.,
$\left(\begin{smallmatrix}X\\ Y\end{smallmatrix}\right)_\varphi$
is a torsionless $\m$-module.

\vskip5pt

The assertion $(\rm ii)$ a direct consequence of the formula $\phi_{\binom{X}{Y}_\varphi} = \widetilde{h}\circ \binom{\phi_X}{\beta^*\phi_Y}$.

\vskip5pt

The assertion $(\rm iii)$ a direct consequence of $(\rm i)$ and $\phi_{\binom{X}{Y}_\varphi} = \widetilde{h}\circ \binom{\phi_X}{\beta^*\phi_Y}$. $\s$



\subsection{Double semi-Gorenstein-projective $\m$-modules} \ \
For a left $\m$-module $\left(\begin{smallmatrix}X\\ Y\end{smallmatrix}\right)_\varphi$
with  $\varphi: Y \longrightarrow X$ a left $A$-map, by Corollary \ref{sgpT2A}(i),
$\left(\begin{smallmatrix}X\\ Y\end{smallmatrix}\right)_\varphi\in \ ^\perp A$ if and only if
the following conditions (1)-(3) hold:

\vskip5pt

(1) \ $X\in \ ^\perp A;$

(2) \ $Y\in \ ^\perp A;$

(3) \ $\varphi^*: X^*\longrightarrow Y^*$ is an epimorphism.

\vskip5pt

By Lemma \ref{T2dualofleftmod}, \ \
$\binom{X}{Y}^{*}_{\varphi} \cong (({\rm Coker} \varphi)^*, \ X^*)_{\pi^*}$ \ \ as right $\m$-modules,
where $\pi^*: ({\rm Coker} \varphi)^*\longrightarrow X^*$ is the right $A$-monomorphism induced by $\pi: X\longrightarrow {\rm Coker}\varphi$. Thus, by the right module version of Corollary \ref{sgpT2A}(i),
\ \ $\left(\begin{smallmatrix}X\\ Y\end{smallmatrix}\right)_\varphi^*\in \ ^\perp A$ if and only if
the following conditions (4)-(6) hold:

\vskip5pt

(4) \ $({\rm Coker} \varphi)^*\in \ ^\perp A;$

(5) \ $X^*\in \ ^\perp A;$

(6) \ $\pi^{**}: X^{**}\longrightarrow ({\rm Coker} \varphi)^{**}$ is an epimorphism.

\begin{lem} \label{doublesgp} \ Let  $\left(\begin{smallmatrix}X\\ Y\end{smallmatrix}\right)_\varphi$
be a left $\m$-module, where $\varphi: Y \longrightarrow X$ is a left $A$-map. Then
$\left(\begin{smallmatrix}X\\ Y\end{smallmatrix}\right)_\varphi$ is double semi-Gorenstein-projective
if and only if the conditions $(1)- (6)$ above hold, and
if and only if the conditions $(1)- (8)$ hold, where

\vskip5pt

$(7)$ \ $Y^*\in \ ^\perp A;$

$(8)$ \ The canonical $A$-map $\beta: {\rm Coker} \pi^*\longrightarrow Y^*$ is an isomorphism.
\end{lem}
\noindent{\bf Proof.}  It remains to show that the conditions (1) - (6) imply the conditions (7) and (8).

\vskip5pt

Assume that the conditions (1) - (6) hold. Since $\varphi^*: X^*\longrightarrow Y^*$ is an epimorphism,
$\beta: {\rm Coker} \pi^*\longrightarrow Y^*$  is an isomorphism  (cf. the diagram (5.2)).
Applying $\Hom_A(-, \ _AA)$ to the exact sequence $0 \longrightarrow ({\rm Coker}\varphi)^* \stackrel {\pi^*} \longrightarrow X^* \stackrel {\varphi^*} \longrightarrow Y^*\longrightarrow 0,$
by the assumption that $\pi^{**}: X^{**}\longrightarrow ({\rm Coker} \varphi)^{**}$ is an epimorphism
and by the assumptions $X^*\in \ ^\perp A$ and $({\rm Coker} \varphi)^*\in \ ^\perp A$, one sees that
$Y^*\in \ ^\perp A.$
$\s$

\vskip5pt

\subsection{A double semi-Gorenstein-projective $\m$-module $M$ with $\phi_M$ monomorphism or epimorphism} \ \

\begin{prop} \label{sgpsandref} \ Let  $\left(\begin{smallmatrix}X\\ Y\end{smallmatrix}\right)_\varphi$
be a left $\m$-module with left $A$-map $\varphi: Y \longrightarrow X$. Then

\vskip5pt

$(\rm i)$ \ \ $\left(\begin{smallmatrix}X\\ Y\end{smallmatrix}\right)_\varphi$ is torsionless and double semi-Gorenstein-projective
if and only if  $\left(\begin{smallmatrix}X\\ Y\end{smallmatrix}\right)_\varphi$ is monic  $($i.e. $\varphi$ is a monomorphism$)$,
\ $X, \ Y$,  and  ${\rm Coker}\varphi$ are double semi-Gorenstein-projective, and $X$ and $Y$ are torsionless.

\vskip5pt

$(\rm ii)$ \ \ $\left(\begin{smallmatrix}X\\ Y\end{smallmatrix}\right)_\varphi$ is double semi-Gorenstein-projective
with epimorphism $\phi_{\left(\begin{smallmatrix}X\\ Y\end{smallmatrix}\right)_\varphi}$ if and only if the following conditions are satisfied$:$

$\bullet$ \ $\varphi^*: X^*\longrightarrow Y^*$ is an epimorphism$;$

$\bullet$ \ All the five modules $X, \ Y,  \ X^*, \ Y^*, \ ({\rm Coker}\varphi)^*$ are semi-Gorenstein-projective$;$

$\bullet$ \ $\phi_X$ and $\phi_Y$ are epimorphisms.

\vskip5pt

$(\rm iii)$ \ \ {\rm(Corollary \ref{sgpT2A})} \ $\left(\begin{smallmatrix}X\\ Y\end{smallmatrix}\right)_\varphi$ is Gorenstein-projective
if and only if $\varphi$ is a monomorphism, $Y$ and ${\rm Coker}\varphi$ are Gorenstein-projective.
\ If this is the case, then $X$ is Gorenstein-projective. \end{prop}
\noindent{\bf Proof.}  $(\rm i)$ \ Assume that $\left(\begin{smallmatrix}X\\ Y\end{smallmatrix}\right)_\varphi$ is torsionless and double semi-Gorenstein-projective.
By Corollary \ref{reflexive}$(\rm i)$,  $\varphi$ is a monomorphism, and $X$ and $Y$ are torsionless.
By Lemma \ref{doublesgp}, all the conditions (1)-(8) hold. Applying $\Hom_A(-, \ _AA)$ to the exact sequence
$0\longrightarrow Y \stackrel \varphi \longrightarrow X \stackrel \pi \longrightarrow {\rm Coker}\varphi \longrightarrow 0$,
since $\varphi^*: X^*\longrightarrow Y^*$ is an epimorphism, and since $X$ and $Y$ are semi-Gorenstein-projective,
it follows that ${\rm Coker}\varphi$ is semi-Gorenstein-projective.

\vskip5pt

Conversely, assume that $\varphi$ is a monomorphism, $X$, \ $Y$ and ${\rm Coker}\varphi$ are double semi-Gorenstein-projective, and that
$X$ and $Y$ are torsionless. By Corollary \ref{reflexive}$(\rm i)$, $\left(\begin{smallmatrix}X\\ Y\end{smallmatrix}\right)_\varphi$ is torsionless.
Again applying $\Hom_A(-, \ _AA)$ to the exact sequence
$0\longrightarrow Y \stackrel \varphi \longrightarrow X \stackrel \pi \longrightarrow {\rm Coker}\varphi \longrightarrow 0$,
since ${\rm Coker}\varphi$ is semi-Gorenstein-projective,  $\varphi^*: X^*\longrightarrow Y^*$ is an epimorphism and
$0 \longrightarrow ({\rm Coker}\varphi)^* \stackrel {\pi^*} \longrightarrow X^* \stackrel {\varphi^*} \longrightarrow Y^*\longrightarrow 0$
is an exact sequence. Since $Y^*$ is semi-Gorenstein-projective, $\pi^{**}: X^{**}\longrightarrow ({\rm Coker} \varphi)^{**}$ is an epimorphism.
Thus, all the conditions (1)-(6) hold. By Lemma \ref{doublesgp}, $\left(\begin{smallmatrix}X\\ Y\end{smallmatrix}\right)_\varphi$ is double semi-Gorenstein-projective.

\vskip10pt

$(\rm ii)$ \ Assume that $\left(\begin{smallmatrix}X\\ Y\end{smallmatrix}\right)_\varphi$ is double semi-Gorenstein-projective
and $\phi_{\left(\begin{smallmatrix}X\\ Y\end{smallmatrix}\right)_\varphi}$ is an epimorphism. By Lemma \ref{doublesgp}, all the conditions (1)-(8) are satisfied.
By Corollary \ref{reflexive}$(\rm ii)$,  $\phi_X$ and $\beta^*\phi_Y$ are epimorphisms,
where $\beta: {\rm Coker}\pi^* \longrightarrow Y^*$ is the canonical $A$-map such that $\varphi^* = \beta p$, \ $p: X^*\longrightarrow {\rm Coker} \pi^*$ is the canonical $A$-epimorphism, and $\beta^*: Y^{**}\longrightarrow ({\rm Coker}\pi^*)^*$ is induced by $\beta$. It remains to show that $\phi_Y$ is an epimorphism.
In fact, by Condition (8),  $\beta^*$ is an isomorphism, hence $\phi_Y$ is an epimorphism.

\vskip5pt

Conversely, assume that $\varphi^*: X^*\longrightarrow Y^*$ is an epimorphism, all the five modules $X$, \ $Y$,  \ $X^*$, \ $Y^*$, \ $({\rm Coker}\varphi)^*$ are semi-Gorenstein-projective, and that
$\phi_X$ and $\phi_Y$ are epimorphisms. Since $\varphi^*$ is an epimorphism, $\beta: {\rm Coker} \pi^*\longrightarrow Y^*$ is an isomorphism (cf. Subsection 5.4),
and hence $\beta^*$ is an isomorphism.  Thus $\beta^*\phi_Y$ is an epimorphism. By Corollary \ref{reflexive}$(\rm ii)$, $\phi_{\left(\begin{smallmatrix}X\\ Y\end{smallmatrix}\right)_\varphi}$ is an epimorphism.

Applying $\Hom_A(-, \ _AA)$ to the exact sequence
$Y \stackrel \varphi \longrightarrow X \stackrel \pi \longrightarrow {\rm Coker}\varphi \longrightarrow 0$,
since $\varphi^*: X^*\longrightarrow Y^*$ is an epimorphism, it follows that
\ $0 \longrightarrow ({\rm Coker}\varphi)^* \stackrel {\pi^*} \longrightarrow X^* \stackrel {\varphi^*} \longrightarrow Y^*\longrightarrow 0$
is an exact sequence. Since $Y^*$ is semi-Gorenstein-projective, $\pi^{**}: X^{**}\longrightarrow ({\rm Coker} \varphi)^{**}$ is an epimorphism.
Thus, all the conditions (1)-(6) hold. By Lemma \ref{doublesgp}, $\left(\begin{smallmatrix}X\\ Y\end{smallmatrix}\right)_\varphi$ is double semi-Gorenstein-projective.

\vskip5pt

$(\rm iii)$ \ This is just Corollary \ref{sgpT2A}. We rewrite here, because in the setting of (i) and (ii),
it admits a simple proof. The ``if" part follows from $(\rm i)$ and $(\rm ii)$ and the fact that Gorenstein-projective modules are closed under extensions.

Assume that $\left(\begin{smallmatrix}X\\ Y\end{smallmatrix}\right)_\varphi$ is Gorenstein-projective. Then by $(\rm i)$ and $(\rm ii)$,
$\varphi$ is a monomorphism, $X$ and $Y$ are Gorenstein-projective, and ${\rm Coker}\varphi$ is double semi-Gorenstein-projective. Moreover, the diagram
$$\xymatrix{0\ar[r] & Y\ar[r]^-{\varphi}\ar[d]^{\phi_Y}_-{\cong} &  X\ar[r]^-{\pi}\ar[d]^{\phi_X}_-{\cong} & {\rm Coker}\varphi \ar[d]^{\phi_{{\rm Coker}\varphi}} \ar[r] & 0  \\
0 \ar[r] & Y^{**}\ar[r]^-{\varphi^{**}} & X^{**} \ar[r]^-{\pi^{**}} & ({\rm Coker}\varphi)^{**} \ar[r] & 0} \eqno (5.3)$$
commutes  with exact rows. So $\phi_{{\rm Coker}\varphi}$ is an isomorphism, and thus ${\rm Coker}\varphi$ is Gorenstein-projective.
$\s$

\subsection{Problems}  \ As remarked in [RZ4, 3.1], all known examples of double semi-Gorenstein-projective modules $M$ such that
$\phi_M$ is a monomorphism (an epimorphism, respectively) are Gorenstein-projective.

\vskip5pt

{\bf Problem 2.} \ Is there a torsionless and double semi-Gorenstein-projective module $M$ such that $M$ is not Gorenstein-projective?

\vskip5pt

{\bf Problem 3.} \ Is there a double semi-Gorenstein-projective module $M$ with  $\phi_M$ an epimorphism  such that $M$ is not semi-Gorenstein-projective?

\vskip5pt

Theorem \ref{sgpsandtorsionless} is a result in this direction.

\subsection{Proof of Theorem \ref{sgpsandtorsionless}}

$(\rm 1)$ \ Assume that any torsionless and double semi-Gorenstein-projective $A$-module is Gorenstein-projective.
Let $M=\left(\begin{smallmatrix}X\\ Y\end{smallmatrix}\right)_\varphi$ be a torsionless and double semi-Gorenstein-projective $\m$-module.
We need to show that $M$ is Gorenstein-projective.

By Proposition \ref{sgpsandref}$(\rm i)$,  $\varphi$ is a monomorphism, $X$, \ $Y$, and ${\rm Coker}\varphi$ are double semi-Gorenstein-projective,
and $X$ and $Y$ are torsionless. By the assumption, $X$ and $Y$ are Gorenstein-projective.

By Lemma \ref{doublesgp}, $\varphi^*$ and $\pi^{**}$ are epimorphisms. Thus, one again has the commutative diagram (5.3) with exact rows,
from which one knows that $\phi_{{\rm Coker}\varphi}$ is also an isomorphism, and hence
${\rm Coker}\varphi$ is Gorenstein-projective.   Thus, $M$ is Gorenstein-projective, by Proposition \ref{sgpsandref}$(\rm iii)$.

\vskip10pt

Conversely, assume that any torsionless and double semi-Gorenstein-projective $\m$-module is Gorenstein-projective.
Let $L$ be a torsionless and double semi-Gorenstein-projective $A$-module. We need to prove that $L$ is Gorenstein-projective.

Since $L$ is torsionless, a left ${\rm add}(A)$-approximation $\varphi: L\longrightarrow P$ of $L$ is a monomorphism, where $P$ is a projective $A$-module.
Since both $P$ and $L$ are semi-Gorenstein-projective and $\varphi$ is a left ${\rm add}(A)$-approximation, it follows that
${\rm Coker}\varphi$ is also semi-Gorenstein-projective. Consider the $\m$-module
$\left(\begin{smallmatrix}P\\ L\end{smallmatrix}\right)_\varphi$. Since $P^*$ and $L^*$ are semi-Gorenstein-projective and
$$0\longrightarrow ({\rm Coker}\varphi)^* \longrightarrow P^* \stackrel {\varphi^*}\longrightarrow L^* \longrightarrow 0$$
is an exact sequence, it follows that $({\rm Coker}\varphi)^*$ is also  semi-Gorenstein-projective. Thus, by Proposition \ref{sgpsandref}$(\rm i)$,
$\left(\begin{smallmatrix}P\\ L\end{smallmatrix}\right)_\varphi$ is a torsionless and double semi-Gorenstein-projective $\m$-module. By the assumption, $\left(\begin{smallmatrix}P\\ L\end{smallmatrix}\right)_\varphi$ is Gorenstein-projective.  Hence  $L$ is Gorenstein-projective, by Proposition \ref{sgpsandref}$(\rm iii)$.

\vskip10pt

$(\rm 2)$ \ \ Assume that any double semi-Gorenstein-projective $A$-module $L$ with $\phi_L$ an epimorphism is Gorenstein-projective.
Let $M=\left(\begin{smallmatrix}X\\ Y\end{smallmatrix}\right)_\varphi$ be a double semi-Gorenstein-projective $\m$-module such that $\phi_M$ is an epimorphism.
We need to prove that $M$ is Gorenstein-projective.

By Proposition \ref{sgpsandref}$(\rm ii)$,  $\varphi^*: X^*\longrightarrow Y^*$ is an epimorphism, all the five modules \ $X$, \ $Y$,  \ $X^*$, \ $Y^*$, \ $({\rm Coker}\varphi)^*$ are semi-Gorenstein-projective,
and $\phi_X$ and $\phi_Y$ are epimorphisms. By the assumption, $X$ and $Y$ are Gorenstein-projective, in particular,
$\phi_X$ and $\phi_Y$ are isomorphisms.

We claim that $\varphi: Y \longrightarrow X$ is a monomorphism and ${\rm Coker}\varphi$ is reflexive. In fact, applying $\Hom_A(-, A)$ to $Y \stackrel \varphi \longrightarrow X \stackrel \pi \longrightarrow {\rm Coker}\varphi \longrightarrow 0$,
since $\varphi^*: X^*\longrightarrow Y^*$ is an epimorphism, it follows that
$$0 \longrightarrow ({\rm Coker}\varphi)^* \stackrel {\pi^*} \longrightarrow X^* \stackrel {\varphi^*} \longrightarrow Y^*\longrightarrow 0$$
is an exact sequence. Since $Y^*$ is semi-Gorenstein-projective,
$$0 \longrightarrow Y^{**} \stackrel {\varphi^{**}} \longrightarrow X^{**} \stackrel {\pi^{**}} \longrightarrow ({\rm Coker}\varphi)^{**}\longrightarrow 0$$
is an exact sequence. Thus, by the functorial property of $\phi$ one has the commutative diagram
$$\xymatrix{& Y\ar[r]^-{\varphi}\ar[d]^{\phi_Y}_{\cong} &  X\ar[r]^-{\pi}\ar[d]^{\phi_X}_{\cong} & {\rm Coker}\varphi \ar[d]^{\phi_{{\rm Coker}\varphi}} \ar[r] & 0  \\
0 \ar[r] & Y^{**}\ar[r]^-{\varphi^{**}} & X^{**} \ar[r]^-{\pi^{**}} & ({\rm Coker}\varphi)^{**} \ar[r] & 0} $$
with exact rows.  Since both $\phi_Y$ and $\varphi^{**}$ are monomorphisms, $\varphi$ is a monomorphism. Also, this commutative diagram shows that $\phi_{{\rm Coker}\varphi}$ is an isomorphism, i.e., ${\rm Coker}\varphi$ is reflexive.
This proves the claim.

Applying $\Hom_A(-, \ _AA)$ to the exact sequence
$0\longrightarrow Y \stackrel \varphi \longrightarrow X \stackrel \pi \longrightarrow {\rm Coker}\varphi \longrightarrow 0$,
since $\varphi^*: X^*\longrightarrow Y^*$ is an epimorphism and $X$ and $Y$ are semi-Gorenstein-projective, it follows that
${\rm Coker}\varphi$ is also semi-Gorenstein-projective. So, ${\rm Coker}\varphi$ is double semi-Gorenstein-projective and reflexive,
i.e.,   ${\rm Coker}\varphi$ is Gorenstein-projective. By Proposition \ref{sgpsandref}$(\rm iii)$,
$M=\left(\begin{smallmatrix}X\\ Y\end{smallmatrix}\right)_\varphi$ is Gorenstein-projective.

\vskip10pt

Conversely, assume that any double semi-Gorenstein-projective $\m$-module $M$ with $\phi_M$ an epimorphism is Gorenstein-projective.
Let $L$ be a double semi-Gorenstein-projective $A$-module such that $\phi_L$ is an epimorphism. We need to show that $L$ is
Gorenstein-projective.

Take a left ${\rm add}(A)$-approximation $\varphi: L\longrightarrow P$ of $L$.
Applying $\Hom_A(-, \ _AA)$ to the exact sequence
$L \stackrel \varphi \longrightarrow P \stackrel \pi \longrightarrow {\rm Coker}\varphi \longrightarrow 0$,
since $\varphi$ is left ${\rm add}(A)$-approximation, $\varphi^*: P^*\longrightarrow L^*$ is an epimorphism and
$$0 \longrightarrow ({\rm Coker}\varphi)^* \stackrel {\pi^*} \longrightarrow P^* \stackrel {\varphi^*} \longrightarrow L^*\longrightarrow 0$$
is an exact sequence. Since $L^*$ and $P^*$ are semi-Gorenstein-projective, so is $({\rm Coker}\varphi)^*$.
Thus, by Proposition \ref{sgpsandref}$(\rm ii)$,
$\left(\begin{smallmatrix}P\\ L\end{smallmatrix}\right)_\varphi$ is a double semi-Gorenstein-projective $\m$-module such that
$\phi_{\left(\begin{smallmatrix}P\\ L\end{smallmatrix}\right)_\varphi}$ is an epimorphism.
By the assumption, $\left(\begin{smallmatrix}P\\ L\end{smallmatrix}\right)_\varphi$ is Gorenstein-projective.  Hence
by Proposition \ref{sgpsandref}$(\rm iii)$, $L$ is Gorenstein-projective. $\s$

\vskip5pt

\section{\bf Double semi-Gorenstein-projective modules which are not monic}

\subsection{Proof of Theorem \ref{doublesgpnotmonic}} \ \ Since by assumption $Y$ is not torsionless and $\varphi: Y\longrightarrow P$ is a left ${\rm add}(A)$-approximation of $Y$,
it follows that $\varphi$ is  not a monomorphism. Thus
$\binom{P}{Y}_\varphi$ is not a monic $T_2(A)$-module. By Corollary \ref{reflexive}(i), $\binom{P}{Y}_\varphi$ is not torsionless.

\vskip5pt

Apply $\Hom_A(-, \ _AA)$ to the exact sequence
\ $Y \stackrel \varphi \longrightarrow P \stackrel \pi \longrightarrow {\rm Coker}\varphi \longrightarrow 0$.
Since $\varphi$ is a left ${\rm add}(A)$-approximation of $Y$, $\varphi^*: P^*\longrightarrow Y^*$ is an epimorphism and
$$0 \longrightarrow ({\rm Coker}\varphi)^* \stackrel {\pi^*} \longrightarrow P^* \stackrel {\varphi^*} \longrightarrow Y^*\longrightarrow 0$$
is an exact sequence. Since $Y^*$ and $P^*$ are semi-Gorenstein-projective, it follows that
$({\rm Coker}\varphi)^*$ is semi-Gorenstein-projective and $\pi^{**}: P^{**}\longrightarrow ({\rm Coker}\varphi)^{**}$ is an epimorphism.
Thus, all the conditions (1) - (6) in Subsection 5.6 are satisfied. By Lemma \ref{doublesgp},
$\left(\begin{smallmatrix}P\\ Y\end{smallmatrix}\right)_\varphi$ is a double semi-Gorenstein-projective $T_2(A)$-module. $\s$

\subsection{A class of double semi-Gorenstein-projective modules which are not torsionless} From now on, $A$ is the algebra $\Lambda(q)$, which has been studied in [RZ2, RZ3], i.e.,
$$A = k\langle x, y, z\rangle / \langle
 x^2,\ y^2,\ z^2,\ yz,\ xy+qyx,\ xz-zx,\ zy-zx\rangle$$
where $q$ is an non-zero element in field $k$, and $q$ is of multiplicative order $\infty$. Then $A$ is a short local algebra of wild representation type,
with a basis $1, \ x, \ y, \ z, \ yx, \ zx$, and with Hilbert type $(|J/J^2|, |J^2|) = (3, 2)$, where $J$ is the Jacobson radical of $A$.
For the studies on short local algebras, we refer to e.g. [L], [Y2], [CV], [AIS],  [RZ5, RZ6].

\vskip5pt

When  $(a,b,c)$ ranges over $k^3\setminus\{0\}$, left $A$-modules $$M(a,b,c): = {}_AA/[A(ax+by+cz) + {\rm soc}A] = A\bar 1$$
give all the $3$-dimensional local $A$-modules. They are $(A/J^2)$-modules.
Since $A/J^2$ is commutative, ${\rm D}(M(a,b,c)) = \Hom_k(M(a,b,c), k)$ are also left $(A/J^2)$-modules, and hence left $A$-modules.
[RZ3, Proposition A.1] asserts that $M(a,b,c)$ and ${\rm D}(M(a,b,c))$ give all the indecomposable left $A$-modules of dimension $3$.

\vskip5pt

For $(a,b,c)\in k^3\setminus\{0\}$, we also consider right $A$-modules $$M'(a,b,c): = A/[(ax+by+cz)A + {\rm soc}A] = \bar 1 A.$$

\begin{lem}\label {RZ2doublesgpnottorsionless} {\rm ([RZ3, 1.7])}  \ An indecomposable $A$-module $M$ of dimension at most $3$ is double semi-Gorenstein-projective which are not torsionless if and only if
$M\cong M(1, -q, c)$ for some $c\in k.$  Moreover, $$M(1, -q, c)^* \cong (x-y)A \cong A/(x-q^{-1}y)A = M'(1, -q^{-1}, 0)$$
where the first isomorphism is given by $f \mapsto f(\bar 1)$, and the second isomorphism is given by $x-y \mapsto \bar 1$.
\end{lem}

\subsection{A class of double semi-Gorenstein-projective $T_2(A)$-modules which are not monic} In order to apply Theorem \ref{doublesgpnotmonic} to get a family of double semi-Gorenstein-projective $T_2(A)$-modules which are not monic, we look for a left ${\rm add}(A)$-approximation of $M(1, -q, c) = A/A(x-qy+cz) = A\bar 1$.  Any $f\in M(1, -q, c)^* = \Hom_A(A\bar 1, \ A)$ is the right multiplication by $f(\bar 1)$. Since ${\rm Im} f\in J$, \ $f(\bar 1) = c_1x+c_2y+c_3z+c_4yx+c_5zx$ with $c_i\in k$, such that $(x-qy+cz)f(\bar 1) = 0$. Thus
$c_1+c_2 = 0, \ c_3 = 0.$ Hence  $f(\bar 1)\in (x-y)A$ and $M(1, -q, c)^*$ has a $k$-basis $f_1, \ f_2, \ f_3$, where
$f_i: M(1, -q, c) \longrightarrow A$  is the left $A$-map given by
$$f_1(\bar 1) = x-y, \ \ \ f_2(\bar 1) = yx, \ \ \ f_3(\bar 1) = zx.$$
Therefore $M(1, -q, c)^* = f_1A$ and
$$f_1: M(1, -q, c) \longrightarrow A$$
is a left ${\rm add}(A)$-approximation of $M(1, -q, c)$. Applying Theorem \ref{doublesgpnotmonic} one gets

\begin{prop} \label{importantexamples} \ For all $c\in k$, the $T_2(A)$-modules $$X(c): = \left(\begin{smallmatrix}_AA\\ M(1, -q, c)\end{smallmatrix}\right)_{f_1}$$
where $f_1: M(1, -q, c) \longrightarrow A$ is the left $A$-map given by $f_1(\bar 1) = x-y,$ are
double semi-Gorenstein-projective, but not monic, and hence not torsionless.

\vskip5pt

Moreover, one has

\vskip5pt

$(\rm i)$ \  $X(c)^* \cong ((x-q^{-1}y)A, \ A_A)_\sigma$, where $\sigma: (x-q^{-1}y)A\longrightarrow A$ is the embedding$;$
$X(c)^{**}\cong \left(\begin{smallmatrix}A \\ A(x-y)A\end{smallmatrix}\right)_{\iota}$ is not semi-Gorenstein-projective, where $\iota: A(x-y)\longrightarrow A$ is the embedding.  In particular,
$X(c)^*$ is not Gorenstein-projective.

\vskip5pt

$(\rm ii)$ \  If one identifies $X(c)^{**}$ with $\left(\begin{smallmatrix}A \\ A(x-y)A\end{smallmatrix}\right)_{\iota}$, then
the canonical $\m$-map $\phi_{X(c)}: X(c)\longrightarrow X(c)^{**}$ reads as
$$\binom{{\rm Id}_A}{r_{x-y}}: \left(\begin{smallmatrix}_AA\\ M(1, -q, c)\end{smallmatrix}\right)_{f_1}\longrightarrow \left(\begin{smallmatrix}_AA\\ A(x-y)A\end{smallmatrix}\right)_{\iota}$$
where $r_{x-y}$ is the right multiplication by $x-y$. Thus, $\phi_{X(c)}$ is neither  a monomorphism nor an epimorphism.
\end{prop}
\noindent{\bf Proof.} \ It remains to prove  $(\rm i)$ and $(\rm ii)$.

\vskip5pt

$(\rm i)$ \ Note that there are isomorphisms
$$({\rm Coker} f_1)^* = (A/A(x-y))^*\cong (x-q^{-1}y)A\cong A/(x-q^{-2}y)A$$ as right $A$-modules, where the first isomorphism is given by $g\mapsto g(\bar1)$, and the second isomorphism is given by $x-q^{-1}y\mapsto \bar 1$.
We stress that $$A/A(x-y)\ncong M(1, -1, 0) = A/[A(x-y)+{\rm soc}A]$$ and that $(A/A(x-y))^* \cong (x-q^{-1}y)A\cong M(1, -1, 0)^*$.
By Lemma \ref{T2dualofleftmod}$(\rm ii)$,  there is a right $\m$-module isomorphism
$$h: \ X(c)^* \cong ((x-q^{-1}y)A, \ A_A)_\sigma$$
where $\sigma: (x-q^{-1}y)A\longrightarrow A$ is the embedding.

\vskip5pt

Note that $({\rm Coker} \sigma)^* = (A/(x-q^{-1}y)A)^* \cong A(x-y)A$ as left $A$-modules, with the isomorphism given by $g\mapsto g(\bar1)$.
By Lemma \ref{doubledual}$(\rm ii)$, there is a left $\m$-module isomorphism
$$\left(\begin{smallmatrix}A \\ A(x-y)A\end{smallmatrix}\right)_{\iota}\cong X(c)^{**}$$
where $\iota: A(x-y)A\longrightarrow A$ is the embedding.

\vskip5pt

Note that $A(x-y)A = A(x-y) \oplus kzx$ is decomposable left $A$-module of dimension $3$.  By [RZ3, Theorem 1.5],  $A(x-y)A$ not semi-Gorenstein-projective. It follows from Corollary \ref{sgpT2A}(1) that
$\left(\begin{smallmatrix}A \\ A(x-y)A\end{smallmatrix}\right)_{\iota}$ is not semi-Gorenstein-projective,
and hence $X(c)^{**}$ is not semi-Gorenstein-projective. In particular,  $X(c)^*$ is not Gorenstein-projective.

\vskip5pt

$(\rm ii)$ \ To get $\phi_{X(c)}$,  we apply Proposition \ref{thecanonicalmap} to $X(c) = \left(\begin{smallmatrix}_AA\\ M(1, -q, c)\end{smallmatrix}\right)_{f_1}$. It is clear that
$\phi_A = {\rm Id}_A$.
Since $X(c)$ is double semi-Gorenstein-projective, it follows from Lemma \ref{doublesgp} that the canonical $A$-map $\beta: {\rm Coker}\sigma \longrightarrow M(1, -q, c)^*$ appeared in Proposition \ref{thecanonicalmap},
is an isomorphism. Without loss of generality, one may regarded $\beta$ as the identity.
Note that $$M(1, -q, c)^{**} \cong M'(1, -q^{-1}, 0)^* \cong A(x-y)A$$
where the first isomorphism is given in Lemma \ref{RZ2doublesgpnottorsionless}, and the second isomorphism is given by
$g \mapsto g(\bar 1)$. By Proposition \ref{thecanonicalmap}(ii), if we identify $X(c)^{**}$ with $\left(\begin{smallmatrix}A \\ A(x-y)A\end{smallmatrix}\right)_{\iota}$, then $\phi_{X(c)}$ reads as
$$\binom{{\rm Id}_A}{\phi_{M(1, -q, c)}}: \left(\begin{smallmatrix}_AA\\ M(1, -q, c)\end{smallmatrix}\right)_{f_1}\longrightarrow \left(\begin{smallmatrix}_AA\\ A(x-y)A\end{smallmatrix}\right)_{\iota}.$$
Since the diagram
$$\xymatrix{
M(1, -q, c)\ar[d]_{\phi_{M(1, -q, c)}} \ar[r]^-{f_1} & A \ar@{=}[d] \\
A(x-y)A\ar[r]^-{\iota} & A}$$
commutes,  it follows that $\phi_{M(1, -q, c)}: M(1, -q, c) \longrightarrow A(x-y)A$ is just the right multiplication $r_{x-y}$.
Thus, $\phi_{M(1, -q, c)}$ is  neither a monomorphism nor an epimorphism, and hence
$\phi_{X(c)}$ is neither a monomorphism nor an epimorphism. $\s$

\vskip20pt

{\bf Acknowledgement}: The author sincerely thanks Claus Michael Ringel for his helpful discussions and comments,
and the anonymous referee for suggestions on the presentation of the paper.

\end{document}